\documentclass[twoside,12pt,a4paper]{article}
\usepackage[marked,extnum,tth,printedin
]{myart98}

\ppnum{}{\htmladdnormallink{math/9807141}%
{http://xxx.lanl.gov/abs/math/9807141}}{1998}

\newcommand{\oper}[1]{\mathcal{#1}}

\titleshort{Wavelets in Banach Spaces}
\authorshort{Vladimir V. Kisil}
\def\tthdump#1{#1}
\tthdump{\typeout{This is MyDef macros file}}
\def\ifundefined#1{\expandafter\ifx\csname#1\endcsname\relax}
\providecommand{\comment}[1]{}

\usepackage{amsfonts}
\tthdump{}
\newcommand{\algebra}[1]{\ensuremath{\mathfrak{#1}}}
\newcommand{\Heisen}[1]{\ensuremath{\mathbb{H}^{#1}}}

\newcommand{\Cliff}[2][\comment]{\ensuremath{%
\mathcal{C}\kern-0.18em\ell(#1,#2)}}
\newcommand{\object}[2][\,]{\ensuremath{\mathrm{#2}#1}}
\newcommand{\Space}[2]{\ensuremath{\mathbb{#1}^{#2}}}
\newcommand{\such}{\,\mid\,}

\newcommand{\FSpace}[2]{\ensuremath{ #1_{#2}}}

\ifundefined{qed}
    \DeclareMathSymbol{\qed}{0}{AMSa}{"03}
\fi
\newcommand{\fourier}[1]{\ensuremath{\mathcal{F}^{#1}}}

\newcommand{\norm}[1]{\left\| #1 \right\|}
\newcommand{\modulus}[1]{\left| #1 \right|}
\newcommand{\scalar}[2]{\left\langle #1,#2 \right\rangle}

\providecommand{\eqref}[1]{\textup{(\ref{#1})}}

\newcommand{\person}[1]{\textsc{#1}}
\newcommand{\mailto}[1]{\href{mailto:#1}{\texttt{#1}}}

\newif\iftth
\providecommand{\href}[2]{#2}

\input{cyr.fd}
\begin{document}
\title{Wavelets in Banach Spaces\thanks{Supported by grant INTAS
93--0322--Ext.}}

\author{\htmladdnormallink{Vladimir V. Kisil}{http://cage.rug.ac.be/~vk/}\\
                           Institute of Mathematics,\\
               Economics and Mechanics,\\
               Odessa State University\\
               ul. Petra Velikogo, 2,\\
               Odessa-57, 270057, UKRAINE\\
               E-mail: \mailto{kisilv@member.ams.org}\\
WWW: \htmladdnormallink{\texttt{http://cage.rug.ac.be/\~{}vk/}}%
{http://cage.rug.ac.be/~vk/}} \date{July 24, 1998} \maketitle
\begin{abstract}
  We describe a construction of wavelets (coherent states) in Banach
  spaces generated by ``admissible'' group representations. Our main
  targets are applications in pure mathematics while connections with
  quantum mechanics are mentioned. As an example we consider operator
  valued Segal-Bargmann type spaces and the Weyl functional calculus.
  \AMSMSC{43A85}{32M99, 43A32, 46E10, 47A60, 47A67, 47C99, 81R30,
    81S10}
\end{abstract}
\newpage {\small \tableofcontents}

\section{Introduction}
\epigraph{Questioning the nature of essence, are we not in danger of
  losing ourselves in the void of the commonplace which suffocates all
  thought?}{M.~Heidegger~\cite{Heidegger61a}.}{}
Wavelets~\cite{BernTayl94,BratJorg97a,Daubechies92,HeilWaln89} (or
coherent
states~\cite{AliAntGazMue,JorgWer94a,Kisil95a,Klauder94b,Klauder94a,%
  KlaSkag85,Perelomov86}) belong to a field of intensive research,
which was carried from many different viewpoints.  This interest is
supported by important applications of the theory in quantum
mechanics, signal processing, etc.  The huge potential of wavelets in
pure mathematics is still to be
explored~\cite{CnopsKisil97a,Kisil95a,Kisil95i,Kisil97a,Kisil97c}.  It
is impossible to give a short overview of all aspect of the theory
here.  We will refer to the recent survey~\cite{AliAntGazMue} which
contains an extensive bibliography (137 entries).

In this paper we investigate a subject which is outside the mainstream
of research: wavelets in Banach space\footnote{In fact, our
  construction is applicable for an arbitrary locally convex
  topological vector space, see Remark~\ref{re:general}.}. The first
paper in this direction known to the present author is the
paper~\cite{Bargmann67} of \person{Bargmann}.  It extends results of
the early paper~\cite{Bargmann61} which describes wavelets in a
Hilbert space associated with a representation of the Heisenberg group
in the Segal-Bargmann space~\cite{Segal60}.  \person{Bargmann}
considered in~\cite{Bargmann67} wavelets which are tempered
distributions.  \person{J.~Klauder} and
\person{J.~McKenna}~\cite{KlauderMcKenna67} extended this idea to
(nontempered) distributions by using a fiducial vector from $\cal D$
rather than from $\cal S$, but they did not publish their results (see
however~\cite{KlaSkag85}).  \person{H.G.~Feichtinger} and
\person{K.H.~Groechenig} in a series of important papers started
by~\cite{FeichGroech89a,FeichGroech89b} described atomic decomposition
in Banach spaces of functions connected with unitary irreducible
representations of groups in Hilbert spaces.  This approach describes
most of the technique used in real function theory and allows to
construct new tools in this area.  Unfortunately the approach has not
received an essential appreciation by function theorists which it
deserves. Note that the Hilbert space technique is the principal tool
for consideration in Banach spaces in the mentioned
papers~\cite{Bargmann67,FeichGroech89a,FeichGroech89b}.

There is also another predecessor from an unexpected field:
combinatorics.  The theory of umbral calculus as it developed by
G.-C.~Rota and his
coauthors~\cite{FoundationIII,KahOdlRota73,RomRota78,Rota64a} contain
all principal formulas for a wavelet transform and its inverse for
particularly selected group $\Space{Z}{}$ (or even convoloid
$\Space{N}{}$~\cite{Kisil97b}).  The key idea of the present paper
(separation a vacuum vector from a test functional in topological or
Banach spaces) is explicitly used in these papers. Due to elementary
nature of this approach its realizations could be probably found in
many other areas of mathematics too.

One could find a reason for which wavelets in Banach space are
underrepresented: two main areas of wavelets application are primary
interested in Hilbert spaces.  In signal processing the finite energy
is defined by the $\FSpace{L}{2}$ norm and in quantum mechanics states
of a systems form a Hilbert space.  However, motivated by pure
mathematics one is interested in Banach or topological spaces also.
One could think on an unified treatment of Hardy spaces
$\FSpace{H}{p}$ without distinguishing value $p=2$ or Segal-Bargmann
type spaces $\FSpace{F}{p}$ (see Subsection~\ref{ss:segbarg-p}).  A
need of wavelets in Banach spaces arises particularly in a study of
principal series representation of semisimple Lie
groups~\cite{Kisil97c}.  Such representations are not square
integrable and therefore the standard theory of wavelets from group
representations~\cite{BernTayl94} is not directly applicable here.
Another appeal arises from operator theory.  To make an advance in
functional calculus of operators we need wavelets defined in Banach
algebras~\cite{Kisil95i,Kisil95e}.  Even from quantum mechanical point
of view Banach spaces are not completely useless: non-unitary
evolution in Banach spaces could describe decay and creation of
particles~\cite{Segal90a}.

In this paper we try to present wavelets in Banach spaces as close as
possible to Hilbert space case but without an explicit use of the
Hilbert space results
(cf.~\cite{Bargmann67,FeichGroech89a,FeichGroech89b}).  An abstract
scheme in Section~\ref{se:banwav} consolidates many results
of~\cite{BernTayl94,FeichGroech89a,FeichGroech89b,KlaSkag85,Perelomov86}
and we receive the standard wavelet theory if assume the space under
consideration to be a Hilbert one.  In Section~\ref{se:algebras} we
apply wavelets for consideration of symbolic and functional calculi of
operators.  As a first application we will consider theory of co- and
contravariant symbols of
operators~\cite{Berezin72,Berezin74,Berezin88} as its direct
consequence.  We give practical examples of wavelets in Banach spaces
in Section~\ref{se:segbarg} by operator valued Segal-Bargmann type
spaces.  Such spaces could be interested for their connections with
the Weyl functional calculus, spectral measures~\cite{Anderson69} and
as a model in quantum mechanics for second quantization of bosonic
fields~\cite{Berezin88,Segal60}. More examples particularly to
wavelets in Banach algebras~\cite{Kisil95e} will be given elsewhere.

We avoid a formulation of the idea in an extreme generality.  Such
formulations while be applicable in more cases require however more
efforts to be applied (and useful!) in any particular situation.  We
select an alternative path: to present the mainstream of the theory
and mention in remarks possible modifications for a generalized
setting (see Subsection~\ref{ss:remarks}).

\section{Coherent States for Banach Spaces}\label{se:coherent}
\label{se:banwav}
\subsection{Abstract Nonsence}\label{ss:abstract}
Let $G$ be a group and $H$ be its closed normal subgroup.  Let $X=G/H$
be the corresponding homogeneous space with an invariant measure
$d\mu$ and $s: X \rightarrow G$ be a Borel section in the principal
bundle $G \rightarrow G/H$.  Let $\pi$ be a continuous representation
of a group $G$ by invertible isometry operators $\pi_g$, $g \in G$ in
a (complex) Banach space $B$.

The following definition simulates ones from the Hilbert space
case~\cite[\S~3.1]{AliAntGazMue}.
\begin{defn} \label{de:coherent1}
  Let $G$, $H$, $X=G/H$, $s: X \rightarrow G$, $\pi: G \rightarrow
  \oper{L}(B)$ be as above. We say that $b_0 \in B$ is a \emph{vacuum
    vector} if for all $h\in H$
  \begin{equation} \label{eq:h-char}
    \pi(h) b_0 = \chi(h) b_0, \qquad \chi(h) \in \Space{C}{}.
  \end{equation}
  We will say that set of vectors $b_x=\pi(x) b_0$, $x\in X$ form a
  family of \emph{coherent states} if there exists a continuous non-zero
  linear functional $l_0 \in B^*$ such that
  \begin{enumerate}
  \item \label{it:norm} $\norm{b_0}=1$, $\norm{l_0}=1$,
    $\scalar{b_0}{l_0}\neq 0$;
  \item \label{it:h-char} $\pi(h)^* l_0=\bar{\chi}(h) l_0$, where
    $\pi(h)^*$ is the adjoint operator to $\pi(h)$;
  \item \label{it:coher-eq} The following equality holds
    \begin{equation} \label{eq:coher-eq}
      \int_X \scalar{\pi(x^{-1}) b_0}{l_0}\, \scalar{\pi(x) b_0}{l_0}\, 
      d\mu(x) = \scalar{b_0}{l_0}.
    \end{equation}
  \end{enumerate}
  The functional $l_0$ is called the \emph{test functional}.  According
  to the strong tradition we call the set $(G,H,\pi,B,b_0,l_0)$
  \emph{admissible} if it satisfies to the above conditions.
\end{defn}
We note that mapping $h \rightarrow \chi(h)$ from~\eqref{eq:h-char}
defines a character of the subgroup $H$.  The following Lemma
demonstrates that condition~\eqref{eq:coher-eq} could be relaxed.
\begin{lem}\label{le:exist1}
  For the existence of a vacuum vector $b_0$ and a test functional
  $l_0$ it is sufficient that there exists a vector $b_0'$ and
  continuous linear functional $l'_0$ satisfying to \eqref{eq:h-char}
  and \textup{\ref{it:h-char}} correspondingly such that the constant
\begin{equation} \label{eq:sq-int}
c = \int_X \scalar{\pi(x^{-1})b'_0}{l'_0}\, \scalar{\pi(x) 
b'_0}{l'_0}\, d\mu(x)
\end{equation}
is non-zero and finite.
\end{lem}
\begin{proof}
  There exist a $x_0\in X$ such that $ \scalar{\pi(x_0^{-1})
    b_0'}{l'_0} \neq 0 $, otherwise one has $c=0$.  Let
  $b_0=\pi(x^{-1}) b_0' \norm{\pi(x^{-1}) b_0'}^{-1} $ and $l_0=l_0'
  \norm{l_0'}^{-1}$. For such $b_0$ and $l_0$ we have \ref{it:norm}
  already fulfilled.  To obtain~\eqref{eq:coher-eq} we change the
  measure $d\mu(x)$.  Let $c_0=\scalar{b_0}{l_0} \neq 0 $ then $d\mu'=
  \norm{\pi(x^{-1}) b_0'} \norm{l_0'} c_0 c^{-1}d\mu$ is the desired
  measure.
\end{proof}
\begin{rem}
  Conditions~\eqref{eq:coher-eq} and~\eqref{eq:sq-int} are known for
  unitary representations in Hilbert spaces as \emph{square
    integrability} (with respect to a subgroup $H$).  Thus our
  definition describes an analog of square integrable representations
  for Banach spaces.  Note that in Hilbert space case $b_0$ and $l_0$
  are often the same function, thus condition~\ref{it:h-char} is
  exactly~\eqref{eq:h-char}.  In the particular but still important
  case of trivial $H=\{e\}$ (and thus $X=G$) all our results take
  simpler forms.
\end{rem}
\begin{conv}
  In that follow we will usually write $x\in X$ and $x^{-1}$ instead
  of $s(x)\in G$ and $s(x)^{-1}$ correspondingly.  The right meaning
  of ``$x$'' could be easily found from the context (whether an
  element of $X$ or $G$ is expected there).
\end{conv}

The wavelet transform (similarly to the Hilbert space case) could be
defined as a mapping from $B$ to a space of bounded continuous
functions over $G$ via representational coefficients
\begin{displaymath}
v \mapsto \widehat{v}(g)= \scalar{\pi(g^{-1})v}{l_0}= 
\scalar{v}{\pi (g)^*l_0}.
\end{displaymath}
Due to~\ref{it:h-char} such functions have simple transformation
properties along orbits $gH$, i.e.
$\widehat{v}(gh)=\bar{\chi}(h)\widehat{v}(g)$, $g\in G$, $h\in H$.
Thus they are completely defined by their values indexed by points of
$X=G/H$.  Therefore we prefer to consider so called reduced wavelet
transform.
\begin{defn}
  The \emph{reduced wavelet transform} $\oper{W}$ from a Banach space
  $B$ to a space of function $\FSpace{F}{}(X)$ on a homogeneous space
  $X=G/H$ defined by a representation $\pi$ of $G$ on $B$, a vacuum
  vector $b_0$ and a test functional $l_0$ is given by the formula
\begin{equation} \label{eq:wave-tr}
\oper{W}: B \rightarrow \FSpace{F}{}(X): v \mapsto \widehat{v}(x)= 
[\oper{W}v] (x)=\scalar{\pi(x^{-1}) v}{l_0}=
\scalar{v}{\pi^*(x)l_0}.
\end{equation}
\end{defn}
There is a natural representation of $G$ in $\FSpace{F}{}(X)$.  For
any $g\in G$ there is a unique decomposition of the form $g=s(x)h$,
$h\in H$, $x\in X$.  We will define $r: G \rightarrow H:
r(g)=h=(s^{-1}(g))^{-1}g$ from the previous equality and write a
formal notation $x=s^{-1}(g)$.  Then there is a geometric action of
$G$ on $X \rightarrow X$ defined as follows
\begin{displaymath}
g: x \mapsto g^{-1} \cdot x = s^{-1} (g^{-1} s(x)).
\end{displaymath}
We define a representation $\lambda(g): \FSpace{F}{}(X) \rightarrow
\FSpace{F}{}(X)$ as follow
\begin{equation} \label{eq:l-rep}
[\lambda(g) f] (x) = \chi(r(g^{-1}\cdot x)) f(g^{-1}\cdot x).
\end{equation}
We recall that $\chi(h)$ is a character of $H$ defined
in~\eqref{eq:h-char} by the vacuum vector $b_0$. For the case of
trivial $H=\{e\}$ \eqref{eq:l-rep} becomes the left regular
representation $\rho_l(g)$ of $G$.
\begin{prop} \label{pr:inter1}
  The reduced wavelet transform $\oper{W}$ intertwines $\pi$ and
  the representation $\lambda$~\eqref{eq:l-rep} on $\FSpace{F}{}(X)$:
\begin{eqnordisp}
  \oper{W} \pi(g) = \lambda(g) \oper{W}.
\end{eqnordisp}
\end{prop}
\begin{proof}
  We have:
\begin{eqnarray*}{}
[\oper{W}( \pi(g) v)] (x) & = & \scalar{\pi(x^{-1}) \pi(g) v }{ l_0} \\
  & = & \scalar{\pi((g^{-1}s(x))^{-1}) v }{ l_0} \\
  & = & \scalar{\pi(r(g^{-1}\cdot x)^{-1})\pi(s(g^{-1}\cdot x)^{-1}) v }{ l_0} \\
  & = & \scalar{\pi(s(g^{-1}\cdot x)^{-1}) v }{\pi^*(r(g^{-1}\cdot x)^{-1}) 
l_0} \\
  & = & \chi(r(g^{-1}\cdot x)^{-1}) [\oper{W} v] (g^{-1}x) \\
  & = & \lambda(g) [\oper{W}v] (x).
\end{eqnarray*}
\end{proof}
\begin{cor}\label{co:pi}
  The function space $\FSpace{F}{}(X)$ is invariant under the
  representation $\lambda$ of $G$.
\end{cor}
We will see that $\FSpace{F}{}(X)$ posses many properties of the
\emph{Hardy space}. 
The duality between $l_0$ and $b_0$ generates a transform dual to $\oper{W}$.
\begin{defn}
  The \emph{inverse wavelet transform} $\oper{M}$ from $
  \FSpace{F}{}(X) $ to $B$ is given by the formula:
\begin{eqnarray}
\oper{M}:  \FSpace{F}{}(X) \rightarrow B: \widehat{v}(x) \mapsto \oper{M} 
[\widehat{v}(x)] & = & \int_X \widehat{v}(x) b_x\,d\mu(x) \nonumber\\
 & = & \int_X \widehat{v}(x) \pi(x)\,d\mu(x) b_0. \label{eq:m-tr}
\end{eqnarray}
\end{defn}
\begin{prop} \label{pr:inter2}
  The inverse wavelet transform $ \oper{M} $ intertwines the
  representation $ \lambda $ on $ \FSpace{F}{}(X) $ and $ \pi $ on
  $B$:
\begin{eqnordisp}
  \oper{M} \lambda(g) = \pi(g) \oper{M}.
\end{eqnordisp}
\end{prop}
\begin{proof}
  We have:
\begin{eqnarray*}
\oper{M} [\lambda(g)\widehat{v}(x)] 
& = & \oper{M} [ \chi(r(g^{-1}\cdot x)) \widehat{v}(g^{-1}\cdot x)] \\
& = & \int_X \chi(r(g^{-1}\cdot x)) \widehat{v}(g^{-1}\cdot x) 
        b_x \,d\mu(x)\\
& = & \chi(r(g^{-1}\cdot x)) \int_X \widehat{v}(x') b_{g\cdot x'}\,d\mu(x')
\\
& = & \pi_g \int_X \widehat{v}(x') b_{x'}\,d\mu(x')\\
& = & \pi_g \oper{M} [\widehat{v}(x')],
\end{eqnarray*}  
where $x'=g^{-1} \cdot x$.
\end{proof}
\begin{cor}\label{co:lambda}
  The image $\oper{M}(\FSpace{F}{}(X))\subset B$ of subspace
  $\FSpace{F}{}(X)$ under the inverse wavelet transform $\oper{M}$ is
  invariant under the representation $\pi$.
\end{cor}
The following proposition explain the usage of the name for
$\oper{M}$.
\begin{thm}
  The operator
\begin{equation} \label{eq:szego}
\oper{P}= \oper{M} \oper{W}: B \rightarrow B
\end{equation}
is a projection of $B$ to its linear subspace for which $b_0$ is
cyclic. Particularly if $\pi$ is an irreducible representation then the
inverse wavelet transform $\oper{M}$ is a \emph{left inverse} operator
on $B$ for the wavelet transform $\oper{W}$:
\begin{eqnordisp}
  \oper{M}\oper{W}=I.
\end{eqnordisp}
\end{thm}
\begin{proof}
  It follows from Propositions~\ref{pr:inter1} and~\ref{pr:inter2}
  that operator $\oper{M}\oper{W}: B \rightarrow B$ intertwines $\pi$
  with itself.  Then Corollaries~\ref{co:pi} and~\ref{co:lambda} imply
  that the image $\oper{M}\oper{W}$ is a $\pi$-invariant subspace of
  $B$ containing $b_0$.  Because $\oper{M}\oper{W}b_0=b_0$ we conclude
  that $\oper{M}\oper{W}$ is a projection.
  
  From irreducibility of $\pi$ by Schur's
  Lemma~\cite[\S~8.2]{Kirillov76} one concludes that
  $\oper{M}\oper{W}=cI$ on $B$ for a constant $c\in\Space{C}{}$.
  Particularly
\begin{displaymath}
\oper{M}\oper{W} b_0= \int_X \scalar{\pi(x^{-1})b_0}{l_0}\, \pi(x) 
b_0\,d\mu(x)=cb_0.
\end{displaymath}
From the condition~\eqref{eq:coher-eq} it follows that
$\scalar{cb_0}{l_0}=\scalar{\oper{M}\oper{W}
  b_0}{l_0}=\scalar{b_0}{l_0}$ and therefore $c=1$.
\end{proof}
We have similar
\begin{thm}
  Operator $\oper{W}\oper{M}$ is a projection of $\FSpace{L}{1}(X)$ to
  $\FSpace{F}{}(X)$. 
\end{thm}

We denote by $\oper{W}^*: \FSpace{F}{}^*(X) \rightarrow B^* $ and
$\oper{M}^*: B^* \rightarrow \FSpace{F}{}^*(X)$ the adjoint (in the
standard sense) operators to $\oper{W}$ and $\oper{M}$ respectively.
\begin{cor}
  We have the following identity:
\begin{equation} \label{eq:isom1}
\scalar{\oper{W} v }{ \oper{M}^* l}_{ \FSpace{F}{}(X) } = \scalar{v}{l}_B, 
\qquad \forall v\in B, \quad l\in B^*
\end{equation}
or equivalently
\begin{equation} \label{eq:isom2}
\int_X \scalar{\pi(x^{-1}) v}{l_0}\, \scalar{\pi(x) b_0}{l}\, d\mu(x) 
= \scalar{v}{l}.
\end{equation}
\end{cor}
\begin{proof}
  We show the equality in the first form~\eqref{eq:isom2} (but will
  apply it often in the second one):
\begin{displaymath}
\scalar{\oper{W} v }{ \oper{M}^* l}_{ \FSpace{F}{}(X) }
 = \scalar{\oper{M}\oper{W} v }{l}_B =\scalar{v}{l}_B.
\end{displaymath}
\end{proof}
\begin{cor}
  The space $ \FSpace{F}{}(X) $ has the reproducing formula
\begin{equation} \label{eq:reprod}
\widehat{v}(y)=\int_X \widehat{v}(x) \, 
   \widehat{b}_0(x^{-1}\cdot y)\,d\mu(x),
\end{equation}
where $ \widehat{b}_0(y)=[\oper{W}b_0] (y)$ is the wavelet transform
of the vacuum vector $b_0$.
\end{cor}
\begin{proof}
  Again we have a simple application of the previous formulas:
\begin{eqnarray}
\widehat{v}(y)& =& \scalar{\pi(y^{-1})v}{l_0} \nonumber \\ 
&= & \int_X \scalar{\pi(x^{-1}) \pi(y^{-1}) v}{l_0}\, \scalar{\pi(x) 
b_0}{l_0}\,d\mu(x) \label{eq:rep-tr1} \\
&= & \int_X \scalar{\pi(s(y\cdot x)^{-1}) v}{l_0}\, \scalar{\pi(x) 
b_0}{l_0}\,d\mu(x) \nonumber \\
&= & \int_X \widehat{v} (y\cdot x)\, \widehat{b}_0(x^{-1})  \,d\mu(x)
\nonumber \\
& = & \int_X \widehat{v}(x)\, \widehat{b}_0(x^{-1}y)\,d\mu(x), \nonumber
\end{eqnarray}
where transformation~\eqref{eq:rep-tr1} is due to~\eqref{eq:isom2}.
\end{proof}
\begin{rem}
  To possess a reproducing kernel---is a well-known property of spaces
  of analytic functions.  The space $\FSpace{F}{}(X)$ shares also
  another important property of analytic functions: it belongs to a
  kernel of a certain first order differential operator with Clifford
  coefficients (the Dirac operator) and a second order operator with
  scalar coefficients (the Laplace
  operator)~\cite{AtiyahSchmid80,Kisil97a,Kisil97c,KnappWallach76}.
\end{rem}
Let us now assume that there are two representations $\pi'$ and
$\pi''$ of the same group $G$ in two different spaces $B'$ and $B''$
such that two admissible sets $(G,H,\pi',B',b_0',l_0')$ and
$(G,H,\pi'',B'',b_0'',l_0'')$ could be constructed for the same normal
subgroup $H\subset G$.
\begin{prop} \label{pr:2repr}
  In the above situation if $F'(X) \subset F''(X)$ then the
  composition $\oper{T}=\oper{M}''\oper{W}'$ of the wavelet transform
  $\oper{W}'$ for $\pi'$ and the inverse wavelet transform
  $\oper{M}''$ for $\pi''$ is an intertwining operator between $\pi'$
  and $\pi''$:
\begin{eqnordisp}[]
  \oper{T}\pi'=\pi''\oper{T}.
\end{eqnordisp}
$\oper{T}$ is defined as follows
\begin{eqnordisp}[eq:T-intertw]
  \oper{T}: b \mapsto \int_X \scalar{\pi'(x^{-1})b}{l'_0}\,
  \pi''(x)b_0''\,d\mu(x).
\end{eqnordisp}
This transformation defines a $B''$-valued linear functional (a
distribution for function spaces) on $B'$.
\end{prop}
The Proposition has an obvious proof.  This simple result is a base
for an alternative approach to functional calculus of
operators~\cite{Kisil95i,Kisil97a} and will be used in
Subsection~\ref{ss:calculus}.  Note also that
formulas~\eqref{eq:wave-tr} and \eqref{eq:m-tr} are particular cases
of~\eqref{eq:T-intertw} because $\oper{W}$ and $\oper{M}$ intertwine
$\pi$ and $\lambda$.

\subsection{Wavelets and a Positive Cone}
The above results are true for wavelets in general.  In applications a
Banach space $B$ is usually equipped with additional structures and
wavelets are interplay with them. We consider an example of such
interaction.

We recall \cite{KreinRutman48}, \cite[Chap.~X]{KantAkil84} the notion
of positivity in Banach spaces. Let $C\subset B$ be a sharp cone, i.e.
$x\in C$ implies that $\lambda x\in C$ and $-\lambda x\not\in C$ for
$\lambda> 0$. We call elements $x\in C$ \emph{positive vectors}, we
say also that $x\geq y$ iff $x-y$ is positive.  There is the
\emph{dual cone} $C^*\in B^*$ defined by the condition
\begin{displaymath}
C^*=\{f\such f\in B^*,\ \scalar{b}{f}\geq 0\ \ \forall x\in C\}.
\end{displaymath}
An operator $A:B\rightarrow B$ is called \emph{positive} if $Ab\geq0$
for all $b\geq0$. If $A$ is positive with respect to $C$ then $A^*$ is
positive with respect to $C^*$.

\begin{defn}
  We call a representation $\pi(g)$ \emph{positive} if there exists a
  vector $b_0\in C$ such that $\pi(x)b_0\in C$ for all $x\in X$.  A
  linear functional $f \in B^*$ is \emph{positive} ($f>0$) with
  respect to a vacuum vector $b_0$ if $\scalar{\pi(x)b_0}{f}\geq 0$
  for all $x\in X$ and $\scalar{\pi(x)b_0}{f}$ is not identically $0$.
\end{defn}
\begin{lem}
  For any positive representation $\pi(g)$ and vacuum vector $b_0$
  there exists a positive test functional.
\end{lem}
\begin{proof}
  Obvious.
\end{proof}
We consider an estimation of positive linear functionals.
\begin{prop} \label{pr:estimation}
  Let $b\in B$ be a vector such that $b=\int_X \widehat{b}(x)
  b_x\,d\mu(x)$.  Let
\begin{enumerate}
\item $D(b)=\{\scalar{\pi(x^{-1})b}{l_0} \such x\in X\}$ be the set of
  value of reduced wavelets transform;
\item $\breve{D}(b)$ be a convex shell of the values of
  $\widehat{b}(x)$;
\item $\hat{D}(b)=\{ \scalar{b}{f} \such f\in C^*, \norm{f}=1, f\geq 0
  \}$.
\end{enumerate}
Then
\begin{eqnordisp}
  D(b) \subset \hat{D}(b) \subset \breve{D}(b).
\end{eqnordisp}
\end{prop}
\begin{proof}
  The first inclusion is obvious. The second could be easily checked:
\begin{displaymath}
\scalar{b}{f}=\scalar{\int_X \widehat{b}(x) b_x\,d\mu(x)}{f}=
\int_X \widehat{b}(x) \scalar{b_x}{f}\,d\mu(x).
\end{displaymath}
\end{proof}

\subsection{Singular Vacuum Vectors} \label{ss:singular}
In many important cases the above general scheme could not be carried
out because the representation $\pi$ of $G$ is not square-integrable
or even not square-integrable modulo a subgroup $H$.  Thereafter the
vacuum vector $b_0$ could not be selected within the original space
$B$ which the representation $\pi$ acts on.  The simplest mathematical
example is the Fourier transform (see Example~\ref{ex:fourier}). In
physics this is the well-known problem of \emph{absence of vacuum
  state} in the constructive algebraic quantum field theory
\cite{Segal90,Segal94,Segal96a}.  The absence of the vacuum within the
linear space of system's states is another illustration to the old
thesis \emph{Natura abhorret vacuum}\footnote{Nature is horrified by
  (any) vacuum (Lat.).} or even more specifically \emph{Natura
  abhorret vectorem vacui}\footnote{Nature is horrified by a carrier
  of nothingness (Lat.).  This illustrates how far a humane beings
  deviated from Nature.}.

We will present a modification of our construction which works in such
a situation.  For a singular vacuum vector the algebraic structure of
group representations could not describe the situation alone and
requires an essential assistance from analytical structures.
\begin{defn} \label{de:coherent2}
  Let $G$, $H$, $X=G/H$, $s: X \rightarrow G$, $\pi: G \rightarrow
  \oper{L}(B)$ be as in Definition~\ref{de:coherent1}. We assume that
  there exist a topological linear space $\widehat{B}\supset B$ such
  that
\begin{enumerate}
\item \label{it:bexist} $B$ is dense in $\widehat{B}$ (in topology of
  $\widehat{B}$) and representation $\pi$ could be uniquely extended
  to the continuous representation $\widehat{\pi}$ on $\widehat{B}$.
\item There exists $b_0 \in \widehat{B}$ be such that for all $h\in H$
\begin{equation} \label{eq:h-char2}
\widehat{\pi}(h) b_0 = \chi(h) b_0, \qquad \chi(h) \in \Space{C}{}.
\end{equation}
\item There exists a continuous non-zero linear functional $l_0 \in
  B^*$ such that \label{it:h-char2a} $\pi(h)^* l_0=\bar{\chi}(h) l_0$,
  where $\pi(h)^*$ is the adjoint operator to $\pi(h)$;
\item \label{it:b2b} The composition $\oper{M}\oper{W}: B \rightarrow
  \widehat{B}$ of the wavelet transform \eqref{eq:wave-tr} and the
  inverse wavelet transform \eqref{eq:m-tr} maps $B$ to $B$.
\item \label{it:coher-eq2} For a vector $p_0\in B$ the following
  equality holds
\begin{equation} \label{eq:coher-eq2}
\scalar{\int_X \scalar{\pi(x^{-1}) p_0}{l_0}\, \pi(x) b_0\, 
d\mu(x) }{l_0}= \scalar{p_0}{l_0},
\end{equation}
where the integral converges in the weak topology of $\widehat{B}$.
\end{enumerate}
As before we call the set of vectors $b_x=\pi(x) b_0$, $x\in X$ by
\emph{coherent states}; the vector $b_0$---a \emph{vacuum vector}; the
functional $l_0$ is called the \emph{test functional} and finally
$p_0$ is the \emph{probe vector}.
\end{defn}
This Definition is more complicated than
Definition~\ref{de:coherent1}.  The equation~\eqref{eq:coher-eq2} is a
substitution for~\eqref{eq:coher-eq} if the linear functional $l_0$ is
not continuous in the topology of $\widehat{B}$.
Example~\ref{ex:fourier} shows that the Definition does not describe
an empty set.  The function theory in \Space{R}{1,1} constructed in
\cite{Kisil97c} provides a more exotic example of a singular vacuum
vector.

We shall show that \ref{it:coher-eq2} could be satisfied by an
adjustment of other components.
\begin{lem}\label{le:exist2}
  For the existence of a vacuum vector $b_0$, a test functional $l_0$,
  and a probe vector $p_0$ it is sufficient that there exists a vector
  $b_0'$ and continuous linear functional $l'_0$ satisfying to
  \textup{\ref{it:bexist}--\ref{it:b2b}} and a vector $p'_0\in B$ such
  that the constant
\begin{eqnordisp}[]
  c=\scalar{\int_X \scalar{\pi(x^{-1}) p_0}{l_0}\, \pi(x) b_0\,
    d\mu(x) }{l_0}
\end{eqnordisp}
is non-zero and finite.
\end{lem}
The proof follows the path for Lemma~\ref{le:exist1}.  The following
Proposition summarizes results which could be obtained in this case.
\begin{prop} Let the wavelet transform $\oper{W}$ \eqref{eq:wave-tr},  
  its inverse $\oper{M}$ \eqref{eq:m-tr}, the representation
  $\lambda(g)$~\eqref{eq:l-rep}, and functional space
  $\FSpace{F}{}(X)$ be adjusted accordingly to
  Definition~\textup{\ref{de:coherent2}}. Then
\begin{enumerate}
\item $\oper{W}$ intertwines $\pi(g)$ and $\lambda(g)$ and the image
  of $\FSpace{F}{}(X)=\oper{W}(B)$ is invariant under $\lambda(g)$.
\item $\oper{M}$ intertwines $\lambda(g)$ and $\widehat{\pi}(g)$ and
  the image of $\oper{M}(\FSpace{F}{}(B))=\oper{M}\oper{W}(B) \subset
  B$ is invariant under $\pi(g)$.
\item If $\oper{M}(\FSpace{F}{}(X))=B$ (particularly if $\pi(g)$ is
  irreducible) then $\oper{M}\oper{W}=I$ otherwise $\oper{M}\oper{W}$
  is a projection $B \rightarrow \oper{M}(\FSpace{F}{}(X))$. In both
  cases $\oper{M}\oper{W}$ is an operator defined by integral
\begin{eqnordisp}[eq:mw2]
  b \mapsto \int_X \scalar{\pi(x^{-1}) b}{l_0}\, \pi(x) b_0\, d\mu(x),
\end{eqnordisp}
\item Space $\FSpace{F}{}(X)$ has a reproducing formula
\begin{eqnordisp}[eq:reprod2]
  \widehat{v}(y)=\scalar{\int_X \widehat{v}(x) \, \pi(x^{-1}y)b_0\,
    dx}{l_0}
\end{eqnordisp}
which could be rewritten as a singular convolution
\begin{eqnordisp}[]
  \widehat{v}(y)=\int_X \widehat{v}(x) \, \widehat{b}(x^{-1}y)\, dx
\end{eqnordisp}
with a distribution $b(y)=\scalar{\pi(y^{-1})b_0}{l_0}$ defined by
\eqref{eq:reprod2}.
\end{enumerate}
\end{prop}
The proof is algebraic and completely similar to
Subsection~\ref{ss:abstract}.

\section{Wavelets in Operator Algebras} \label{se:algebras}
We are going to apply the above abstract scheme to special spaces,
which our main targets---wavelets on operator algebras. This gives a
possibility to study operators by means of functions---symbols of
operators.

\subsection{Co- and Contravariant Symbols of Operators}
We construct a realization of the wavelet transform as co- and
contravariant symbols (also known as Wick and anti-Wick symbols) of
operators.  These symbols and their connections with wavelets in
Hilbert spaces are known for a
while~\cite{Berezin72,Berezin74,Berezin75,Berezin88}.  However their
realization (described bellow) as wavelets in Banach algebras seems to
be new.

Let $\pi(g)$ be a representation of a group $G$ in a Banach space $B$
by isometry operators. Then we could define two new representations
for groups $G$ and $G \times G$ correspondingly in the space
$\oper{L}(B)$ of bounded linear operators $B\rightarrow B$:
\begin{eqnarray}
\widehat{\pi}&:& G \rightarrow \oper{L}(\oper{L}(B)) : 
 A \mapsto \pi(g)^{-1} A \pi(g), \label{eq:repG}\\
\widetilde{\pi}&:& G \times G \rightarrow \oper{L}(\oper{L}(B)) :
 A \mapsto \pi(g_1)^{-1} A \pi(g_2), 
\label{eq:repGG}
\end{eqnarray}
where $A\in \oper{L}(B)$. Note that $\widehat{\pi}(g)$ are algebra
automorphisms of $\oper{L}(B)$ for all $g$. Representation
$\widetilde{\pi}(g_1,g_2)$ is an algebra homomorphism from
$\oper{L}(B)$ to the algebra $\FSpace{L}{(g_1,g_2)}(B)$ of linear
operators on $B$ equipped with a composition
\begin{displaymath}
A_1\circ A_2= A_1\,\pi(g_1)^{-1}\pi(g_2)\,A_2
\end{displaymath}
with the usual multiplication of operators in the right-hand side. The
r\^ole of such algebra homomorphisms in a symbolical calculus of
operators was explained in~\cite{Howe80b}. It is also obvious that
$\widehat{\pi}(g)$ is the restriction of $\widetilde{\pi}(g_1,g_2)$ to
the diagonal of $G \times G$.

Let there are selected a vacuum vector $b_0\in B$ and a test
functional $l_0\in B^*$ for $\pi$.  Then there are the canonically
associated vacuum vector $P_0 \in \oper{L}(B)$ and test functional
$f_0 \in\FSpace{L}{}^*(B)$ defined as follows:
\begin{eqnarray}
P_0 &:& B \rightarrow B: b \mapsto P_0b=\scalar{b}{l_0}b_0;
\label{eq:def-P0}\\
f_0 &:& \oper{L}(B) \rightarrow \Space{C}{}: A \mapsto 
\scalar{Ab_0}{l_0}. \label{eq:def-f0}
\end{eqnarray}
They define the following coherent states and transformations of the
test functional
\begin{eqnarray*}
P_g&=&\widehat{\pi}(g)P_0=\scalar{\cdot}{l_g}b_g, \qquad
P_{(g_1,g_2)}=\widetilde{\pi}(g_1,g_2)P_0=\scalar{\cdot}{l_{g_1}}b_{g_2}, 
\\
f_g&=&\widehat{\pi}^*(g)f_0=\scalar{\,\cdot\, b_g}{l_g}, \qquad
f_{(g_1,g_2)}=\widehat{\pi}^*(g_1,g_2)f_0=\scalar{\,\cdot\, 
b_{g_1}}{l_{g_2}}
\end{eqnarray*}
where as usually we denote $b_g=\pi(g)b_0$, $l_g=\pi^*(g)l_0$. All
these formulas take simpler forms for Hilbert spaces if $l_0=b_0$.
\begin{defn}
  The \emph{covariant (pre-)symbol} $A(x)$ ($A(x_1,x_2)$) of an
  operator $A$ acting on a Banach space $B$ defined by $b_0\in B$ and
  $l_0 \in B^*$ is its wavelet transform with respect to
  representation $\widehat{\pi}(g)$~\eqref{eq:repG}
  ($\widetilde{\pi}(g_1,g_2)$~\eqref{eq:repGG}) and the functional
  $f_0$~\eqref{eq:def-f0}, i.e.  they are defined by the formulas
\begin{eqnarray}
\hspace*{-5mm} A(x)&=& (\widehat{\pi}(x)A,f_0)= \scalar{\pi(x)^{-1}A\pi(x) 
b_0}{l_0} =\scalar{A b_x}{l_x}, \label{eq:covariant}\\
\hspace*{-5mm} A(x_1,x_2)&=& (\widetilde{\pi}(x_1,x_2)A,f_0)= 
\scalar{\pi(x_1)^{-1}A\pi(x_2)b_0}{l_0} 
=\scalar{A b_{x_2}}{l_{x_1}}.
\end{eqnarray}
The \emph{contravariant (pre-)symbol} of an operator $A$ is a function
$\breve{A}(x)$ (a function $\breve{A}(x_1,x_2)$ correspondingly) such
that $A$ is the inverse wavelet transform of $\breve{A}(x)$ (of
$\breve{A}(x_1,x_2)$ correspondingly) with respect to
$\widehat{\pi}(g)$ ($\widetilde{\pi}(g_1,g_2)$), i.e.
\begin{eqnarray}
A&=&\int_X \breve{A}(x) \widehat{\pi}(x) P_0\,d\mu(x)
=\int_X \breve{A}(x) P_x\,d\mu(x),\\
A&=&\int_X \int_X \breve{A}(x_1,x_2) \widetilde{\pi}(x_1,x_2) 
P_0\,d\mu(x_1)\,d\mu(x_2) \nonumber\\
&=&\int_X \int_X \breve{A}(x_1,x_2) P_{(x_1,x_2)}\,d\mu(x_1)\,d\mu(x_2),
\end{eqnarray}
where the integral is defined in the weak sense.
\end{defn}
Obviously the covariant symbol $\breve{A}(x)$ is the restriction of
the covariant pre\-sym\-bol $\breve{A}(x_1,x_2)$ to the diagonal of $G
\times G$.
\begin{prop}
  Mapping $\sigma: A \mapsto \sigma_A(x_1,x_2)$ of operators to their
  covariant symbols is the algebra homomorphism from algebra of
  operators on $B$ to algebra of integral operators on $
  \FSpace{F}{}(G)$, i.e.
\begin{equation} \label{eq:symbol}
\sigma_{A_1A_2}(x_1,x_3)=\int_X \sigma_{A_1}(x_1,x_2) 
\sigma_{B_2}(x_2,x_3)\,
d\mu(x_2).
\end{equation}
\end{prop}
\begin{proof}
  One could easily see that:
\begin{eqnarray}
\lefteqn{\int_X \sigma_{A_1}(x_1,x_2) \sigma_{A_2}(x_2,x_3)\, d\mu(x_2)} 
\qquad && \nonumber\\
 & = & \int_X \scalar{\pi(x_1)A_1\pi(x_2^{-1}) b_0}{l_0}\,
\scalar{\pi(x_2)A_2\pi(x_3^{-1}) b_0}{l_0}\,d\mu(x_2) \nonumber \\
& = & \int_X \scalar{\pi(x_2^{-1}) b_0}{A_1^*\pi^*(x_1)l_0}\,
\scalar{\pi(x_2)A_2\pi(x_3^{-1}) b_0}{l_0}\,d\mu(x_2) \nonumber \\
& = & \scalar{A_2\pi(x_3^{-1}) b_0}{A_1^*\pi^*(x_1)l_0} \label{eq:tr1-1} 
\\
& = & \scalar{\pi(x_1)A_1A_2\pi(x_3^{-1}) b_0}{l_0} \nonumber \\
& = &  \sigma_{A_1A_2}(x_1,x_3), \nonumber
\end{eqnarray}
where transformation~\eqref{eq:tr1-1} is due to~\eqref{eq:isom2}.
\end{proof}
The following proposition is obvious.
\begin{prop}
  An operator $A$ could be reconstructed from its covariant presymbol
  $A(g_1,g_2)$ by the formula
\begin{eqnordisp}
  Av=\int _G \int_G A(g_1,g_2) \widehat{v}(g_2)\,d\mu(g_2)
  b_{g_1}\,d\mu(g_1).
\end{eqnordisp}
\end{prop}
We have a particular interest in operators closely connected with the
representation $\pi_g$.
\begin{prop}
  Let an operator $A$ on $B$ is defined by the formula
\begin{eqnordisp}
  Av=\int_G a(g)\, \pi_g v\, d\mu(g)
\end{eqnordisp}
for a function $a(g)$ on $G$. Then $A' \oper{W} = \oper{W} A$ where
$A'$ is a two-sided convolution on $G$ defined by the formula
\begin{eqnordisp}{}
  [A' \widehat{v}] (h) = \int_G \int_G a(g_1) \widehat{b}_0(g_2)
  \widehat{v}(g_1^{-1} h g_2)\, d\mu(g_1)\, d\mu(g_2).
\end{eqnordisp}
\end{prop}
For operator algebras there are the standard notions of positivity:
any operator of the form $A^*A$ is positive; if algebra is realized as
operators on a Hilbert space $H$ then $b\in H$ defines a positive
functional $f_b(A)=\scalar{Ab}{b}$.  Thus the following proposition is
a direct consequence of the Proposition~\ref{pr:estimation}.
\begin{prop} \textup{\cite[Thm.~1]{Berezin72}}
  Let $A$ be an operator, let $D(A)$ be the set of values of the
  covariant symbol $A(x)$, let $\breve{D}(A)$ be a convex shell of the
  values of contravariant symbol $\breve{A}(x)$.  Let $\hat{D}(A)$ be
  the set of values of the quadratic form $\scalar{Ab}{b}$ for all
  vectors $\norm{b}=1$.  Then
\begin{eqnordisp}
  D(b) \subset \hat{D}(b) \subset \breve{D}(b).
\end{eqnordisp}
\end{prop}
\begin{example}
  There are at least two very important realizations of symbolical
  calculus of operators. The theory of pseudodifferential operators
  (PDO)~\cite{Cordes79,Shubin87,MTaylor81} is based on the
  Schr\"odinger representation of the Heisenberg group $\Space{H}{n}$
  (see Subsection~\ref{ss:schrodinger}) on the spaces of functions
  $\FSpace{L}{p}(\Space{R}{n})$~\cite{Howe80b}. The Wick and anti-Wick
  symbolical calculi~\cite{Berezin72,Berezin88} arise from the
  Segal-Bargmann representation~\cite{Segal60,Bargmann61} (see
  Subsection~\ref{ex:seg-barg}) of the same group $\Space{H}{n}$.
  Connections (intertwining operators) between these two
  representations were exploited in~\cite{Howe80b} to obtain
  fundamentals of the theory of PDO.
\end{example}

\subsection{Functional Calculus and Group Representations}
\label{ss:calculus}
This Subsection illustrates a new approach to functional calculus of
operators outlined in~\cite{Kisil95i,Kisil97a}. The approach uses the
intertwining property for two representations instead of an algebraic
homomorphism.

Let $\algebra{B}$ be a Banach algebra and $\mathbf{T}\subset
\algebra{B}$ be its subset of elements.  Let $G$ be a group, $H$ be
its normal subgroup and $X=G/H$---the corresponding homogeneous space.
We assume that there is a representation $\tau$ depending from
$\mathbf{T}\subset \algebra{B}$ defined on measurable functions from
$\FSpace{L}{}(X,\algebra{B})$ by the formula
\begin{eqnordisp}[eq:tau-def]
  \tau(g) f(x)= t(g,x) f(g^{-1}\cdot x), \qquad f(x)\in
  \FSpace{L}{}(X,\algebra{B}),
\end{eqnordisp}
where $t(g,x): \algebra{B} \rightarrow \algebra{B}$ depends from $x\in
X$ and $g\in G$.  It is convenient to use a linear functional $l \in
\algebra{B}'$ to make the situation more tractable by reducing it to
the scalar case. Using $l$ we could define a representation
$\tau_l(x)$ on $\FSpace{F}{}(X)$ by the following formula
\begin{eqnordisp}[eq:tau-scal]
  \tau_l(x): f_l(y)=\scalar{f(y)}{l} \mapsto
  [\tau_l(x)f_l](y)=\scalar{\tau(x) f(y)}{l},
\end{eqnordisp}
where $ f(y) \in \FSpace{L}{}(X,\algebra{B})$, $l\in \algebra{B}'$.
We will understand convergence of all integrals involving $\tau$ in a
weak sense, i.e.  as convergence of all corresponding integrals with
$\tau_l$, $l\in \algebra{B}'$.  We also say that $\tau$ is
\emph{irreducible} if all $\tau_l$ are irreducible.
\begin{rem}
  If $\algebra{B}$ is realized as an algebra of operators on a Banach
  space $B$ then $l\in\algebra{B}'$ could be realized as an element of
  $B\otimes B'$.  In this case the formula~\eqref{eq:tau-scal} looks
  like~\eqref{eq:covariant}.  The important difference is the
  following.  In~\eqref{eq:covariant} the representation in the
  operator algebra $\algebra{B}$ arises from a representation in
  Banach space $B$ and is the same for all elements of $\algebra{B}$.
  Representation $\tau_l$ in~\eqref{eq:tau-scal} is defined via the
  representation $\tau$ which depends in its turn from a set
  $\mathbf{T}\subset \algebra{B}$.  Such representations are usually
  connected with some (non-linear) geometric actions of a group
  directly on operator algebra.  Examples of these geometric actions
  are the representation of the Heisenberg group~\eqref{eq:rep-a}
  leading to the Weyl functional calculus, fractional-linear
  transformations of operators leading~\cite{Kisil97a} to
  Dunford-Riesz functional calculus and monogenic functional
  calculus~\cite{Kisil95i}.  Thus such representations contain
  important information on $\mathbf{T}$.
\end{rem}

We also assume that there is a representation $\pi$ of $G$ in
$\FSpace{F}{}(\Omega)$ with a vacuum vector $b_0$, a test functional
$l_0$ and the system of wavelets (coherent states) $b_x$, $x\in X$,
which were main actors in the previous Section.  Let
$\pi^*(g)=\pi(g^{-1})^*$ be the adjoint representation of $\pi(g)$ in
$\FSpace{F}{}'(\Omega)$.

We need a preselected element $T_0(x) \in \FSpace{L}{}(X,\algebra{B})$
which plays a r\^ole of a \emph{vacum vector} for the representation
$\tau$, it is defined by the condition:
\begin{eqnordisp}[eq:t-vac-def]
  \int_X \widehat{b}_0 (x')\, \tau(x') T_0(x)\,dx'=T_0(x),
\end{eqnordisp}
where $\widehat{b}_0(x)=\scalar{\pi(x^{-1}) b_0}{l_0}$ is the wavelet
transform of the vacuum vector $b_0\in \FSpace{F}{}(\Omega)$ for
$\pi$.

\begin{lem}
  A vacuum vector for $\tau$ always exists and is given by the formula
\begin{eqnordisp}[eq:t-vac-exp]
  T_0(x)=\int_X \widehat{b}_0(x') \, \tau(x') T(x)\, dx',
\end{eqnordisp}
where $T(x)\in \FSpace{L}{}(X,\algebra{B})$ is an arbitrary element
which the integral~\eqref{eq:t-vac-exp} converges for.  If $\tau$ is
irreducible (i.e.  all $\tau_l$~\eqref{eq:tau-scal} are irreducible)
in the linear span of $\tau(x')T(x)$, $x'\in X$ then $T_0(x)$ does not
depend from a particular chose of $T(x)$.
\end{lem}
\begin{proof}
  First we could easily verify condition~\eqref{eq:t-vac-def} for the
  $T_0$ defined by~\eqref{eq:t-vac-exp}:
\begin{eqnarray}
\lefteqn{ \int_X \widehat{b}_0 (x)\, \tau(x) T_0(x_1)\,dx 
 =  \int_X \widehat{b}_0 (x)\, \tau(x) \int_X \widehat{b}_0(x') \, 
\tau(x') T(x_1)\, dx'\,dx }\hspace{3cm}&& \nonumber \\
& = & \int_X \int_X  \widehat{b}_0 (x)\, \widehat{b}_0(x') \, 
 \tau(x)\tau(x') T(x_1)\, dx'\,dx \nonumber\\
& = & \int_X \left( \int_X  \widehat{b}_0 (x) \widehat{b}_0(x^{-1}x'') 
\,dx\right) \tau(x'') T(x_1)\, dx''\quad \label{eq:t-vac1}\\
& = & \int_X  \widehat{b}_0 (x'')\, \tau(x'') T(x_1)\, dx'' 
\label{eq:t-vac2}\\
& = & T_0(x_1). \nonumber
\end{eqnarray}
Here we use the change of variables $x''=x\cdot x'$
in~\eqref{eq:t-vac1} and reproducing property~\eqref{eq:reprod} of
$\widehat{b}_0(x)$ in~\eqref{eq:t-vac2}.

To prove that for any admissible $T(x)$ we will receive the same
$T_0(x)$ is enough to pass from the representation $\tau$ to
representations $\tau_l$~\eqref{eq:tau-scal} defined by $l\in
\algebra{B}'$.  Then we deal with scalar valued (not operator valued)
functions and knew that one could use any admissible vector
$T_l(x)=\scalar{T(x)}{l}$ as a vacuum vector in the reconstruction
formula~\eqref{eq:m-tr}.
\end{proof}

Now we could specify the Definition 1.1 from~\cite{Kisil95i} as
follows.
\begin{defn}\label{de:calc}
  Let $G$, $H$, $X=G/H$, $\algebra{B}$, $\mathbf{T}$, $\tau$, $\pi$,
  $\FSpace{F}{}$, $b_0$, $T_0$ be as described above.  One says that a
  continuous linear $\algebra{B}$-valued functional
  $\Phi_{\mathbf{T}}(\cdot,x): \FSpace{F}{} (\Omega) \rightarrow
  \algebra{B}$, parametrized by a point $x\in X$ and depending from
  $\textbf{T}\subset \algebra{B}$:
\begin{eqnordisp}[]
  \Phi_{\mathbf{T}}(\cdot,x): f(y) \mapsto [\Phi_{\mathbf{T}}f](x)=
  \int_\Omega f(y) \Phi_{\mathbf{T}} (y,x)\,dy
\end{eqnordisp}
is a \emph{functional calculus} if
\begin{enumerate}
\item $\Phi_{\mathbf{T}}$ is an intertwining operator between $\pi(g)$
  and $\tau(g)$, namely
\begin{eqnordisp}[eq:phi-inter1]
  [\Phi_{\mathbf{T}}\pi(g)f(y)](x)= \tau(g)[\Phi_{\mathbf{T}}
  f(y)](x),
\end{eqnordisp}
for all $g\in G$ and $f(y)\in \FSpace{F}{}(\Omega)$
\item $\Phi_{\mathbf{T}}$ maps the vacuum vector $b_0(y)$ for the
  representation $\pi$ to the vacuum vector $T_0(x)$ for the
  representation $\tau$:
\begin{eqnordisp}[eq:phi-vac]
  [\Phi_{\mathbf{T}}b_0(y)] (x)=T_0(x).
\end{eqnordisp}
\end{enumerate}
$\algebra{B}$-valued distribution $\Phi_{\mathbf{T}}(y, x_0)$,
$s(x_0)=e\in G$ associated with $\algebra{B}$-valued linear functional
on $\FSpace{F}{}(\Omega)$ is called a \emph{spectral decomposition} of
operators $\mathbf{T}$.
\end{defn}

Representation $\tau$ in~\eqref{eq:phi-inter1} is defined
by~\eqref{eq:tau-def}:
\begin{eqnordisp}[]
  \tau(g) [\Phi_{\mathbf{T}} f(y)]( x)= t(g,x) [\Phi_{\mathbf{T}}
  f(y)](g^{-1}\cdot x).
\end{eqnordisp}
We could state~\eqref{eq:phi-inter1} equivalently as
\begin{eqnordisp}[]
  [I_y \otimes \tau_x(g)] \Phi(y,x) = [\pi_y^*(g^{-1}) \otimes I_x]
  \Phi(y,x).
\end{eqnordisp}
\begin{rem}\label{re:zero}
  The functional calculus $\Phi_{\mathbf{T}}(y,x)$ as defined here has
  the explicit covariant property with respect to variable $x$.  Thus
  it could be restored by the representation $\tau$ from a single
  value, e.g.  $\Phi_{\mathbf{T}}(y, s^{-1}(e))$, where $e$ is the
  identity of $G$.  We particularly will calculate only
  $[\Phi_{\mathbf{T}}f](s^{-1}(e))$ in Subsection~\ref{ss:oper} as the
  value of a functional calculus.  This value is usually denoted by
  $f(\mathbf{T})$ and is exactly the functional calculus of operators
  in the traditional meaning.
\end{rem}
In particular cases different characteristics of the spectral
decomposition could give relevant information on the set of operators
$\mathbf{T}$, e.g.  the support $\object{supp}\!_y
\Phi_{\mathbf{T}}(y,x_0)$ of $\Phi_{\mathbf{T}}(y,x_0)$
\begin{eqnordisp}[]
  f(y)=0\ \forall y\in \object{supp}\!_y\Phi_{\mathbf{T}}(y,x_0) \quad
  \Rightarrow\quad [\Phi_{\mathbf{T}}f](x_0)=0
\end{eqnordisp}
is called \emph{(joint) spectrum} of set $\mathbf{T}\subset
\algebra{B}$.  This definition of the spectrum is connected with the
\emph{Arveson-Connes spectral
  theory}~\cite{Arveson80,Connes73,Takesaki83} while there are several
important differences mentioned in~\cite[Rem.~4.4]{Kisil95i}.

In the paper~\cite{Kisil95i} the approach was illustrated by a newly
developed functional calculus for several non-commuting operators
based on M\"obius transformations of the unit ball in $\Space{R}{n}$.
It was shown in~\cite[\S~7]{Kisil97a} that the classic Dunford-Riesz
functional calculus is generated by a representation of
$SL(2,\Space{R}{})$ within this procedure.  However an abstract scheme
of the approach was not presented yet.  We give some its elements
here.

From Proposition~\ref{pr:2repr} we know a general form of an
intertwining operator of two related representations of a group, which
could be employed here.  Let $l_0(y)$ be the distribution
corresponding to a test functional $l_0$ for the representation $\pi$
on $\FSpace{F}{}(\Omega)$ such that we could write
\begin{eqnordisp}[]
  \scalar{f(y)}{l_0}_{\FSpace{F}{}(\Omega)} = \int_\Omega f(y)
  l_0(y)\,dy.
\end{eqnordisp}
We denote also by $\pi^*(x) l_0(y)$, $x\in X$, $y\in \Omega$
distributions corresponding to linear functionals $\pi^*(x) l_0$,
where $\pi^*(x)$ is the adjoint representation to $\pi$ on the space
$\FSpace{F}{}(\Omega)$.

\begin{prop}[Spectral syntesis]\label{pr:calc}
  Under assumption of Proposition \ref{pr:2repr} the functional
  calculus exists and is unique.  The spectral decomposition
  $\Phi_{\mathbf{T}}(y,x)$ as a distribution on $\Omega$ is given by
  the formula
\begin{eqnordisp}[eq:calc-dis]
  \Phi_{\mathbf{T}}(y,x)=\int_X \pi^*(x) l_0(y) \, \tau(x) T_0(x)\,
  dx.
\end{eqnordisp}
The functional calculus $\Phi_{\mathbf{T}}(\cdot,x)$ as a mapping
$\FSpace{F}{}(\Omega) \rightarrow \algebra{B}$ is given
correspondingly
\begin{eqnordisp}[eq:calc-map]
  \Phi_{\mathbf{T}}(\cdot,x): f(y) \mapsto [\Phi_{\mathbf{T}}f(y)](x)=
  \int_X \int_\Omega \scalar{f(y)}{\pi^*(x') l_0(y)}\, \tau(x') T_0(x)\,dy\,
  dx'.
\end{eqnordisp}
\end{prop}
\begin{proof}
  Obviously \eqref{eq:calc-dis} and~\eqref{eq:calc-map} are
  equivalent. Thus we will prove \eqref{eq:calc-map} only.  For an
  arbitrary $f(y)\in \FSpace{F}{}(\Omega)$ we could write
\begin{eqnarray}
[\Phi_{\mathbf{T}} f(y)](x) & = & \left[\Phi_{\mathbf{T}} 
  \int_X \widehat{f}(x') \pi(x') b_0(y)\,dx' \right](x)
  \label{eq:phi-ded1}\\
& = & \int_X \widehat{f}(x') [\Phi_{\mathbf{T}} \pi(x') b_0(y)](x)\,dx'
  \label{eq:phi-ded2}\\
& = & \int_X \widehat{f}(x') \tau(x') [\Phi_{\mathbf{T}} b_0(y)](x)\,dx'
  \label{eq:phi-ded3}\\
& = & \int_X \widehat{f}(x')\, \tau(x') T_0(x)\,dx' \label{eq:phi-ded4}\\
& = & \int_X \scalar{f(y)}{\pi^*(x') l_0}\, \tau(x') T_0(x)\,dx' 
\label{eq:phi-ded5}
\end{eqnarray}
We use in~\eqref{eq:phi-ded1} that functions in $\FSpace{F}{}(\Omega)$
are superpositions of coherent states,
transformation~\eqref{eq:phi-ded2} is made by linearity and continuity
of $\Phi_{\mathbf{T}}$, step~\eqref{eq:phi-ded3} is due to
condition~\eqref{eq:phi-inter1} and we finally apply
\eqref{eq:phi-vac} to receive \eqref{eq:phi-ded4}.  Thus it is proven
that the functional calculus which is continuous, linear, and
satisfies to~\eqref{eq:phi-inter1} and~\eqref{eq:phi-vac} (if exists)
is unique and given by~\eqref{eq:phi-ded5}.  Now we should check
that~\eqref{eq:phi-ded5} really gives the right answer.

We will check first that~\eqref{eq:phi-ded5} satisfies to
\eqref{eq:phi-inter1}:
\begin{eqnarray}
\tau(g) \Phi(y,x) & =& \tau(x_1) \int_X \pi^*(x') l_0(y)\, \tau(x') 
T_0(x)\,dx' \nonumber\\
& =& \int_X \pi^*(x') l_0(y)\,  \tau(g\cdot x') T_0(x)\,dx' \nonumber\\
& =& \int_X \pi^*(g^{-1}\cdot x'') l_0(y)\,  \tau(x'') T_0(x)\,dx'' 
\label{eq:phi-ded6}\\
& =& \pi^*(g^{-1}) \int_X \pi^*(x'') l_0(y)\,  \tau(x'') T_0(x)\,dx'' 
\nonumber\\
& =& \pi^*(g^{-1}) \Phi(y,x) \nonumber,
\end{eqnarray}
where we made substitution $x''=g\cdot x'$ in~\eqref{eq:phi-ded6}.
Finally \eqref{eq:phi-vac} directly follows from the
condition~\eqref{eq:t-vac-def}.
\end{proof}
Let there exists $L_0(x)\in \FSpace{L}{}'(X,\algebra{B})$---a test
functional for a vacuum vector $T_0(x)$ and representation $\tau$,
i.e.
\begin{eqnordisp}[]
  \scalar{T_0}{L_0}_{\FSpace{L}{}(X,\algebra{B})} = \int_X
  \scalar{\tau(x^{-1} T_0)}{L_0}_{\FSpace{L}{}(X,\algebra{B})}
  \scalar{\tau(x) T_0}{L_0}_{\FSpace{L}{}(X,\algebra{B})}\,dx,
\end{eqnordisp}
where
\begin{eqnordisp}[]
  \scalar{T_0}{L_0}_{\FSpace{L}{}(X,\algebra{B})} = \int_X
  \scalar{T_0(x)}{L_0(x)}_{\algebra{B}}\,dx
\end{eqnordisp}
and $\scalar{T_0(x)}{L_0(x)}_{\algebra{B}}$ is the pairing between
$\algebra{B}$ and $\algebra{B}'$.
\begin{prop}[Spectral analysis]
  If a $\algebra{B}$-valued function $F(x)$ from
  $\FSpace{L}{}(X,\algebra{B})$ belongs to the closer of the linear
  span of $\tau(x')T_0(x)$, $x'\in X$ then
\begin{eqnordisp}[]
  F(x)= [\Phi_{\mathbf{T}} f(y)] (x),
\end{eqnordisp}
where
\begin{eqnordisp}[eq:sp-anal]
  f(y) = \int_X
  \scalar{\tau(x^{-1}F)}{L_0}_{\FSpace{L}{}(X,\algebra{B})}\,
  \pi(x)b_0(y)\,dx.
\end{eqnordisp}
\end{prop}
\begin{proof}
  The formula~\eqref{eq:sp-anal} is just another realization of
  intertwining operator~\eqref{eq:T-intertw} from
  Proposition~\ref{pr:2repr}.
\end{proof}

Let $K: \FSpace{F}{}(\Omega_2) \rightarrow \FSpace{F}{}(\Omega_2)$ be
an intertwining mapping between two representations $\pi_1$ and
$\pi_2$ of groups $G_1$ and $G_2$ in spaces $\FSpace{F}{}(\Omega_1)$
and $\FSpace{F}{}(\Omega_1)$ respectively. Let $K(z,y)$, $z\in
\Omega_1$, $y\in \Omega_2$ be the Schwartz kernel of $K$.
\begin{thm}[Mapping of spectral decompositions]\label{th:sp-map}
  Let
\begin{eqnarray*}
f_1(z)= [Kf_2](z) & = & \int_{\Omega_2} f_2(y) K(z,y)\,dy\\
\Phi_{\mathbf{T}_1} (z,x) & = & \int_{\Omega_2} K(z,y) 
\Phi_{\mathbf{T}_2}(y,x)\,dy,
\end{eqnarray*}
where a functional calculus $\Phi_{\mathbf{T}_2}$ is defined by
representations $\pi_2$ and $\tau$. Then $\Phi_{\mathbf{T}_2}(z,x)$ is
a functional calculus for $\pi_1$ and $\tau$ and we have an identity:
\begin{eqnordisp}[eq:sp-map]
  [\Phi_{\mathbf{T}_1} f_1(z)](x)= [\Phi_{\mathbf{T}_2} f_2(y)] (x).
\end{eqnordisp}
\end{thm}
\begin{proof}
  The intertwining property for $\Phi_{\mathbf{T}_2}(z,x)$ follows
  from transitivity.  The identity~\eqref{eq:sp-map} is a simple
  application of the Fubini theorem:
\begin{eqnarray}
[\Phi_{\mathbf{S}} g(z)](x) &=& \int_{\Omega_1} g(z) \Phi_{\mathbf{S}} 
(z,x)\,dz \nonumber\\
&=& \int_{\Omega_1} g(z) \int_{\Omega_2} K(z,y) \Phi_{\mathbf{T}}(y,x) 
\,dy\,dz \nonumber\\
&=& \int_{\Omega_2} \int_{\Omega_1} g(z) K(z,y) \, dz\,
\Phi_{\mathbf{T}}(y,x) \,dy \nonumber\\
&=& \int_{\Omega_2} g(f(y)) \Phi_{\mathbf{T}}(y,x) \,dy \nonumber\\
&=& [\Phi_{\mathbf{T}} g(f(y))](x). \nonumber
\end{eqnarray}
\end{proof}
This Theorem could be turned in the spectral mapping theorem under
suitable conditions~\cite[Thm.~3.19]{Kisil95i}.

\subsection{Remarks on Future Generalizations}\label{ss:remarks}
We would like to make some final remarks and mention possible
generalizations.
\begin{rem}\label{re:general}
  Our consideration could be applied even for a more general setting:
  coherent states in an arbitrary locally convex topological vector
  space.  The transition to topological case is straightforward:
  anywhere one should substitute norm estimations by topological
  conditions in the standard way.
\end{rem}
\begin{rem}
  Irreducibility of $\pi$ guarantees that a vacuum vector and a test
  functional exist, even any vector could be taken as a vacuum and
  coherent states will be dense in whole $B$. But irreducibility is
  not necessary for the theory in fact.
\end{rem}
\begin{example}
  One could take a normal base $\{ b_j\}$, $j\in \Space{Z}{}$ in $B$
  indexed by integers (both positive and negative) and take a
  biorthogonal system of linear functionals $\{ l_j\}$:
  $\scalar{b_j}{l_k}=\delta_{jk}$.  Then the group \Space{Z}{} of
  integers has a representation $\sigma$ on $B$, which is defined on
  the base as follows $\sigma(k)b_j=b_{j+k}$, $k\in \Space{Z}{}$.  The
  adjoint representation is $\sigma^*(k)l_j=l_{j-k}$, $k\in
  \Space{Z}{}$.  Then $b_0$ and $l_0$ could be taken as a vacuum
  vector and test functional correspondingly and $b_0$ is cyclic for
  $B$.  All results of this Subsection~\ref{se:coherent} are true
  then.  However the representation $\sigma$ is far from irreducible
  one.  One could find a proof that cyclic vectors exist for a wide
  family of reducible representations in~\cite{CotlSad90b}.
\end{example}
\begin{rem}
  Many of our results (e.g.,~\eqref{eq:symbol}) are based only on the
  formula~\eqref{eq:coher-eq} and are not related to a group structure
  within the set $X=G/H$, which indexed coherent states.  Thus one
  could give a general definition of coherent states based only
  on~\eqref{eq:coher-eq} interpreted as a \emph{resolution of unity}.
  For Hilbert space such general definitions are used
  in~\cite{Berezin72,Berezin88,Klauder94b,Klauder95a}.
  
  However other results like~\eqref{eq:reprod} use the group structure
  and future incites could be given by the harmonic analysis of
  particular groups involved into subject. For this reasons we stay
  with our Definition~\ref{de:coherent1}, which is not being of the
  extreme generality could give more information in particular cases.
\end{rem}
\begin{rem} \label{re:linearization}
  We often meet in applications non-linear representations.  But such
  representations could be linearized in the standard way.  Let a
  non-linear representation $\pi$ of $G$ be defined in a linear space
  $B$ such that $\pi(g): B \rightarrow B$ are uniformly bounded for
  all $g\in G$.  In such a case for any $G$-homogeneous space $X$ with
  an invariant measure $d\mu(x)$ all linear integrals of the form

\begin{eqnordisp}[eq:integrals]
  f(b)=\int_X f(x) \pi(x) b\,d\mu(x), \qquad f\in\FSpace{L}{1}(X),\ 
  b\in B, \norm{b}=1
\end{eqnordisp}
converge in $B$.  They form an invariant linear subspace $B_X$ of $B$
with respect to an associated linear representation $\widetilde{\pi}$
of $G$ (or even of the convolution algebra $\FSpace{L}{1}(G)$) defined
as follows:
\begin{eqnordisp}[eq:rep-op-ext]
  \widetilde{\pi}(g): B_X \rightarrow B_X : f(b) \mapsto
  \widetilde{\pi}(g) f(b)=\int_X f(g\cdot x) \pi(x) b\,d\mu(x),
\end{eqnordisp}
where $f\in\FSpace{L}{1}(G),\ b\in B$.  To make
integrals~\eqref{eq:integrals} convergent for locally convex
topological spaces one may restrict symbols $f(x)$ to smooth functions
with compact support.  Then many of previous results could be
translated to this setting.
\end{rem}
\begin{rem}
  There is a powerful way to generalize the above construction. One
  could define the wavelet transform $\oper{W}: b \rightarrow
  \scalar{\pi(x^{-1})b}{l_0}$ not via a scalar \Space{C}{}-valued
  product with a linear functional $l_0$ but as $L_0(\pi(x^{-1})b)$
  for a linear map $L_0: B \rightarrow \algebra{B}$ where
  $\algebra{B}$ is the space of bounded linear operators on $B$. Then
  the wavelet transform $\widehat{b}(x)$ is $\algebra{B}$-valued
  function on $X$. The inverse wavelet transform then could be again
  defined by the formula~\eqref{eq:m-tr}. For $\algebra{B}$-valued
  wavelet transform $\chi(h)$ should not be a $\Space{C}{}$-valued
  character of $H$. Instead of it one could consider arbitrary
  representation of $H$ in $B$ such that $L_0(\pi(h)b)=\chi(h)L_0(b)$
  for all $b\in B$, $h\in H$. As a simple example of such a wavelet
  transform one could consider $\Space{R}{1,1}$-valued function theory
  constructed in~\cite{Kisil97c}.
\end{rem}
\begin{rem}
  Another generalization is possible within a framework of
  multiresolution wavelet theory~\cite{BratJorg97a}, which is very
  useful for certain reducible representation~\cite{CnopsKisil97a}.
  Namely, one considers not a single vacuum vector $b_0$ but a family
  of such vectors $b_{0,j}$, $1\leq j \leq n$ and associated set of
  test functionals $l_{0,j}$.  Coherent states $b_{x,j}$ are indexed
  now by $x\in X$ and $1\leq j \leq n$.  We have a family of wavelet
  transforms $\oper{W}_j$:
\begin{eqnordisp}{}
  \widehat{b}(j,x)=[\oper{W}_j b](x) = \scalar{\pi(x^{-1})
    b}{l_{0,j}}, \qquad x\in X,\ 1\leq j \leq n
\end{eqnordisp}
and single reconstruction formula
\begin{eqnordisp}[]
  \oper{M}\widehat{b}(j,x) = \int_X \sum_{j=1}^n\widehat{b}(j,x)
  \pi(x) b_{0,j}\,dx.
\end{eqnordisp}
It is interesting to note that this reconstruction formula could be
found already in papers of \person{M.G.~Krein}~\cite{Krein48a}.
\end{rem}

\section{Examples}
\label{se:segbarg}
We are going to demonstrate that the above construction is not only
algebraically attractive but also belongs to the heart of analysis.
More examples could be found
in~\cite{CnopsKisil97a,Kisil97a,Kisil97c} and will be given elsewhere.

\subsection{The Heisenberg Group and Schr\"odinger Representation} 
\label{ss:schrodinger}

We will consider a realization of the previous results in a particular
cases of the Fourier transform and
Segal-Bargmann~\cite{Bargmann61,Segal60} type spaces
$\FSpace{F}{p}(\Space{C}{n})$.  They arise from representations of the
Heisenberg group $\Space{H}{n}$~\cite{Folland89,Howe80a,MTaylor86} on
$\FSpace{L}{p}(\Space{R}{n})$.

The Lie algebra $\algebra{h}_n$ of $\Space{H}{n}$ spanned by $\{T,
P_j, Q_j\}$, $n=1,\ldots,n$ is defined by the commutation relations:
\begin{equation}\label{eq:a-heisenberg}
[P_i,Q_j]=T\delta_{ij}.
\end{equation}
They are known from quantum mechanics as the canonical commutation
relations of coordinates and momentum operators.  An element $g\in
\Heisen{n}$ could be represented as $g=(t,z)$ with $t\in\Space{R}{}$,
$z=(z_1,\ldots,z_n)\in \Space{C}{n}$ and the group law is given by
\begin{equation}\label{eq:g-heisenberg}
g*g'=(t,z)*(t',z')=(t+t'+\frac{1}{2}\sum_{j=1}^n\Im(\bar{z}_j 
z_j'), z+z'),
\end{equation}
where $\Im z$ denotes the imaginary part of a complex number $z$.  The
Heisenberg group is (non-commutative) nilpotent step 2 Lie group.

We take a representation of $\Space{H}{n}$ in
$\FSpace{L}{p}(\Space{R}{n})$, $1<p<\infty$ by operators of shift and
multiplication~\cite[\S~1.1]{MTaylor86}:
\begin{equation} \label{eq:schrodinger}
g=(t,z): f(y) \rightarrow [\sigma_{(t,z)}f](y)=e^{i(2t-\sqrt{2}vy+uv)} 
f(y- \sqrt{2}u), \qquad z=u+iv,
\end{equation}
i.e., this is the Schr\"odinger type representation with parameter
$\hbar=1$. These operators are isometries in
$\FSpace{L}{p}(\Space{R}{n})$ and the adjoint representation
$\pi^*_{(t,z)}=\pi_{(-t,-z)}$ in $\FSpace{L}{q}(\Space{R}{n})$,
$p^{-1}+q^{-1}=1$ is given by a formula similar to
\eqref{eq:schrodinger}.
 
\subsection{Wavelet Transforms for the Heisenberg Group in Function Spaces} 
\label{ss:segbarg-p}
\begin{example} \label{ex:fourier}
  We start from the subgroup $H=\Space{R}{n+1}=\{(t,z)\such
  \Im(z)=0\}$.  Then $X=G/H=\Space{R}{n}$ and an invariant measure
  coincides with the Lebesgue measure.  Mappings $s: \Space{R}{n}
  \rightarrow \Space{H}{n}$ and $r: \Space{H}{n} \rightarrow H$ are
  defined by the identities $s(x)=(0,ix)$, $s^{-1}(t,z)=\Im z$,
  $r(t,u+iv)=(t,u)$.  The composition law $s^{-1}((t,z)\cdot
  s(x))=x+u$ reduces to Euclidean shifts on $\Space{R}{n}$.  We also
  find $s^{-1}((s(x_1))^{-1}\cdot s(x_2))=x_2-x_1$ and
  $r((s(x_1))^{-1}\cdot s(x_2))= 0$.
  
  We consider the representation $\sigma(g)$ of \Space{H}{n} in the
  space of smooth rapidly decreasing functions
  $B=\mathcal{S}(\Space{R}{n})$.  As a character of $H=\Space{R}{n+1}$
  we take the $\chi(t,u)=e^{2it}$.  The corresponding test functional
  $l_0$ satisfying to \ref{it:h-char2a} is the integration
  $l_0(f)=(2\pi)^{-n/2}\int_{\Space{R}{n}} f(y)\,dy$.  Thus the
  wavelet transform is as follows
\begin{eqnordisp}[eq:fourier]
  \widehat{f}(x)=\int_{\Space{R}{n}} \sigma(s(x)^{-1}) f(y)\,dy =
  (2\pi)^{-n/2}\int_{\Space{R}{n}} e^{i \sqrt{2}xy} f(y)\,dy
\end{eqnordisp}
and is nothing else but the Fourier transform\footnote{The
  \emph{inverse} Fourier transform in fact. In our case the signs
  selection is opposite to the standard one, but we will neglect this
  difference.}.

Now we arrive to the absence of a vacuum vector in $B$, indeed there
is no a $f(x)\in \mathcal{S}(\Space{R}{n})$ such that
\begin{displaymath}
[\sigma(t,u) f](y)= \chi(t,u)f(y) \iff 
e^{i2t} f(y- \sqrt{2}u)=e^{i2t} f(y).
\end{displaymath}
There is a way out accordingly to Subsection~\ref{ss:singular}.  We
take $B'=\FSpace{L}{\infty}(\Space{R}{n})\supset B$ and the vacuum
vector $b_0(y)\equiv (2\pi)^{-n/2} \in B'$. Then coherent states are
$b_x(y)=(2\pi)^{-n/2} e^{-i \sqrt{2}xy}$ and the inverse wavelet
transform is defined by the inverse Fourier transform
\begin{displaymath}
f(y)= \int_{\Space{R}{n}} \widehat{f}(y) b_x(y)\,dx
= (2\pi)^{-n/2} \int_{\Space{R}{n}} \widehat{f}(y) e^{-i \sqrt{2}xy}\,dx.
\end{displaymath}
The condition \ref{it:b2b} $\oper{M}\oper{W}: B \rightarrow B$ follows
from the composition of two facts $\oper{W}: B\rightarrow B$ and
almost identical to it $\oper{M}: B\rightarrow B$, which are proved in
standard analysis textbooks (see for
example~\cite[\S~IV.2.3]{KirGvi82}).  To check scaling
\eqref{eq:coher-eq2} according to the tradition in analysis
~\cite{Howe80a} we take a probe vector $p_0=e^{-y^2/2}\in B$.  Due to
well known formula $\int_{-\infty}^{+\infty}e^{-y^2/2}dy=(2\pi)^{1/2}$
of real analysis we have
\begin{eqnarray*}
{\scalar{\int_X \scalar{\widehat{p}_0(x)}{l_0}b_x\,dx}{l_0}}
&=& (2\pi)^{-n} \int\!\int\!\int 
e^{i\sqrt{2}xy} e^{-y^2/2}\,dy e^{-i\sqrt{2}xw}\,dx\,dw\\
&=&(2\pi)^{n/2} \int_{\Space{R}{n}} e^{-y^2/2}dy\\
&=& \scalar{p_0}{l_0}.
\end{eqnarray*}
Thus our scaling is correct.  $\oper{W}$ and $\oper{M}$ intertwine the
left regular representation --- multiplication by $e^{i \sqrt[]{2}yv}$
with operators
\begin{eqnarray*}                          
[\lambda(g) f] (x) &=& \chi(r(g^{-1}\cdot x)) f(g^{-1}\cdot x)\\
&=& e^{i \sqrt[]{2}\cdot 0} f(x-\sqrt[]{2}u)=f(x-\sqrt[]{2}u),
\end{eqnarray*}
i.e. with Euclidean shifts.  From the identity $\scalar{\oper{W} v }{
  \oper{M}^* l}_{ \FSpace{F}{}(X) } =
\scalar{v}{l}_B$~\eqref{eq:isom1} follows the Plancherel's identity:
\begin{eqnarray*}
\int_{\Space{R}{n}} \widehat{v}(y) \widehat{l}(y)\, dy &=&
\int_{\Space{R}{n}} {v}(x) {l}(x)\, dx .
\end{eqnarray*}
These are basic and important properties of the Fourier transform.

The Schr\"odinger representation is irreducible on
$\mathcal{S}(\Space{R}{n})$ thus $\oper{M}=\oper{W}^{-1}$.  Thereafter
integral formulas \eqref{eq:mw2} and \eqref{eq:reprod2} representing
operators $\oper{M}\oper{W}=\oper{W}\oper{M}=1$ correspondingly give
an integral resolution for a convolution with the Dirac delta
$\delta(x)$.  We have integral resolution for the Dirac delta
\begin{eqnordisp}[]
  \delta(x-y)=(2\pi)^{-n/2}\int_{\Space{R}{n}} e^{i\xi(x-y)}\,d\xi.
\end{eqnordisp}

All described results on the Fourier transform are a part of any
graduate curriculum.  What is a reason for a reinvention of a bicycle
here?  First, the same path works with minor modifications for a
function theory in \Space{R}{1,1} described in~\cite{Kisil97c}.
Second we will use this interpretation of the Fourier transform in
Example~\ref{ex:weyl} for a demonstration how the Weyl functional
calculus fits in the scheme outlined in~\cite{Kisil95i,Kisil97a} and
Subsection~\ref{ss:calculus}.
\begin{rem}\label{re:ax+b}
  Of course, the Heisenberg group is not the only possible source for
  the Fourier transform.  We could consider the ``$ax+b$''
  group~\cite[Chap.~7]{MTaylor86} of the affine transformations of
  Euclidean space $\Space{R}{n}$.  The normal subgroup $H=\Space{R}{}$
  of dilations generates the homogeneous space $X=\Space{R}{n}$ on
  which shifts act simply transitively.  The Fourier transform deduced
  from this setting will naturally exhibit scaling properties.  We
  could alternatively consider a group $\Space{M}{n}$ of M\"obius
  transformation~\cite[Chap.~2]{Cnops94a} in $\Space{R}{n+1}$ which
  map upper half plane to itself.  Then there is an induced action of
  $\Space{M}{n}$ on $\Space{R}{n}$---the boundary of upper half plane.
  $\Space{M}{n}$ generated by composition of the affine
  transformations and the Kelvin inverse~\cite[Chap.~2]{Cnops94a}.  If
  we take the normal subgroup $H$ generated by dilations and the
  Kelvin inverse then the quotient space $X$ will again coincide with
  $\Space{R}{n}$ and we immediately arrive to the above case.  On the
  other hand the Fourier transform derived in such a way could be
  easily connected with the plane wave decomposition~\cite{Sommen88}
  in Clifford analysis~\cite{BraDelSom82,DelSomSou92}.
\end{rem}
\end{example}
\begin{example}\label{ex:seg-barg}
  As a subgroup $H$ we select now the center of $\Space{H}{n}$
  consisting of elements $(t,0)$.  Of course $\Omega=G/H$ isomorphic
  to $\Space{C}{n}$ and mapping $s: \Space{C}{n} \rightarrow G$ simply
  is defined as $s(z)=(0,z)$.  The Haar measure on $ \Space{H}{n} $
  coincides with the standard Lebesgue measure on
  $\Space{R}{2n+1}$~\cite[\S~1.1]{MTaylor86} thus the invariant
  measure on $\Omega$ also coincides with the Lebesgue measure on
  $\Space{C}{n}$.  Note also that composition law $s^{-1}(g\cdot
  s(z))$ reduces to Euclidean shifts on $\Space{C}{n}$.  We also find
  $s^{-1}((s(z_1))^{-1}\cdot s(z_2))=z_2-z_1$ and
  $r((s(z_1))^{-1}\cdot s(z_2))= \frac{1}{2} \Im \bar{z}_1z_2$.
  
  As a ``vacuum vector'' we will select the original \emph{vacuum
    vector} of quantum mechanics---the Gauss function
  $f_0(x)=e^{-x^2/2}$ which belongs to all
  $\FSpace{L}{p}(\Space{R}{n})$.  Its transformations are defined as
  follow:
\begin{eqnarray*}
f_g(x)=[\pi_{(t,z)} f_0](x) & = & e^{i(2t-\sqrt{2}vx+uv)}\,e^{-{(x- 
\sqrt{2}u)}^2/2}\\
     & = & e^{2it-(u^2+v^2)/2} e^{- ((u-iv)^2+x^2)/2+\sqrt{2}(u-iv)x} 
\\
     & = & e^{2it-z\bar{z}/2}e^{- (\bar{z}^2+x^2)/2+\sqrt{2}\bar{z}x}.
\end{eqnarray*}
Particularly $[\pi_{(t,0)} f_0](x)=e^{-2it}f_0(x)$, i.e., it really is
a vacuum vector in the sense of our definition with respect to $H$.
For the same reasons we could take $l_0(x)=e^{-x^2/2} \in
\FSpace{L}{q}(\Space{R}{n})$, $p^{-1}+q^{-1}=1$ as the test
functional.

It could be shown that $[\pi_{(0,z)} f_0](x)$ belongs to
$\FSpace{L}{q}(\Space{R}{n})\otimes \FSpace{L}{p}(\Space{C}{n})$ for
all $p>1$ and $q>1$, $p^{-1}+q^{-1}=1$.  Thus
transformation~\eqref{eq:wave-tr} with the kernel $[\pi_{(0,z)}
f_0](x)$ is an embedding $\FSpace{L}{p}(\Space{R}{n}) \rightarrow
\FSpace{L}{p}(\Space{C}{n})$ and is given by the formula
\begin{eqnarray}
\widehat{f}(z)&=&\scalar{f}{\pi_{s(z)}f_0}\nonumber \\
     &=&\pi^{-n/4}\int_{\Space{R}{n}} f(x)\, e^{-z\bar{z}/2}\,e^{-
(z^2+x^2)/2+\sqrt{2}zx}\,dx \nonumber \\
     &=&e^{-z\bar{z}/2}\pi^{-n/4}\int_{\Space{R}{n}} f(x)\,e^{-
(z^2+x^2)/2+\sqrt{2}zx}\,dx \label{eq:tr-bargmann}.
\end{eqnarray}
Then $\widehat{f}(g)$ belongs to $\FSpace{L}{p}( \Space{C}{n} , dg)$
or its preferably to say that function
$\breve{f}(z)=e^{z\bar{z}/2}\widehat{f}(t_0,z)$ belongs to space
$\FSpace{L}{p}( \Space{C}{n} , e^{- \modulus{z}^2 }dg)$ because
$\breve{f}(z)$ is analytic in $z$.  Such functions for $p=2$ form the
\emph{Segal-Bargmann space} $ \FSpace{F}{2}( \Space{C}{n}, e^{-
  \modulus{z}^2 }dg) $ of functions~\cite{Bargmann61,Segal60}, which
are analytic by $z$ and square-integrable with respect the Gaussian
measure $e^{- \modulus{z}^2}dz$.  For this reason we call the image of
the transformation \eqref{eq:tr-bargmann} by \emph{Segal-Bargmann type
  space} $\FSpace{F}{p}( \Space{C}{n}, e^{- \modulus{z}^2 }dg)$.
Analyticity of $\breve{f}(z)$ is equivalent to condition $( \frac{
  \partial }{ \partial\bar{z}_j } + \frac{1}{2} z_j I )
\widehat{f}(z)=0 $.

The integral in~\eqref{eq:tr-bargmann} is the well-known
Segal-Bargmann transform~\cite{Bargmann61,Segal60}. Inverse to it is
given by a realization of~\eqref{eq:m-tr}:
\begin{eqnarray}
f(x) & = & \int_{ \Space{C}{n} } \widehat{f}(z) f_{s(z)}(x)\,dz 
\nonumber\\
& = & \int_{ \Space{C}{n} } \widehat{f}(u,v) e^{iv(u-\sqrt{2}x)}\,e^{-{(x-
\sqrt{2}u)}^2/2} \,du\,dv \label{eq:sb-inverse}\\
& = & \int_{
\Space{C}{n} } \breve{f}(z) e^{- (\bar{z}^2+x^2)/2+\sqrt{2}\bar{z}x}\, e^{-
\modulus{z}^2}\, dz.  \nonumber
\end{eqnarray}
The corresponding operator $\oper{P}$~\eqref{eq:szego} is an identity
operator $ \FSpace{L}{p}(\Space{R}{n}) \rightarrow
\FSpace{L}{p}(\Space{R}{n}) $ and~\eqref{eq:szego} gives an integral
presentation of the Dirac delta.

Integral transformations~\eqref{eq:tr-bargmann}
and~\eqref{eq:sb-inverse} intertwines the Schr\"odinger
representation~\eqref{eq:schrodinger} with the following realization
of representation~\eqref{eq:l-rep}:
\begin{eqnarray}
\lambda(t,z) f(w) & = & \widehat{f}_0(z^{-1}\cdot w) 
\bar{\chi}(t+r(z^{-1}\cdot w))\\
& = & \widehat{f}_0(w-z)e^{it+i\Im(\bar{z}w)} 
\end{eqnarray}

Meanwhile the orthoprojection $\FSpace{L}{2}( \Space{C}{n}, e^{-
  \modulus{z}^2 }dg) \rightarrow \FSpace{F}{2}( \Space{C}{n}, e^{-
  \modulus{z}^2 }dg)$ is of a separate interest and is a principal
ingredient in Berezin quantization~\cite{Berezin88,Coburn94a}.  We
could easy find its kernel from~\eqref{eq:reprod}.  Indeed,
$\widehat{f}_0(z)=e^{ - \modulus{z}^2 }$, then the kernel is
\begin{eqnarray*}
K(z,w) & = & \widehat{f}_0(z^{-1}\cdot w) \bar{\chi}(r(z^{-1}\cdot w))\\
& = & \widehat{f}_0(w-z)e^{i\Im(\bar{z}w)} \\
& = & \exp\left(\frac{1}{2}(- \modulus{w-z}^2 +w\bar{z}-z\bar{w})\right)\\
& = &\exp\left(\frac{1}{2}(- \modulus{z}^2- \modulus{w}^2) 
+w\bar{z}\right).
\end{eqnarray*}
To receive the reproducing kernel for functions
$\breve{f}(z)=e^{\modulus{z}^2} \widehat{f}(z) $ in the Segal-Bargmann
space we should multiply $K(z,w)$ by $e^{(-\modulus{z}^2+
  \modulus{w}^2)/2}$ which gives the standard reproducing kernel $=
\exp(- \modulus{z}^2 +w\bar{z})$ \cite[(1.10)]{Bargmann61}.
\end{example}

\subsection{Operator Valued Representations of the Heisenberg Group}\label{ss:oper}
We proceed now with our main targets: wavelets in operator algebras.
We shell show that well-known and new functional calculi are
realizations of the scheme from Subsection~\ref{ss:calculus}.
\begin{conv} \textup{\cite{Anderson69}} \label{co:oper}
  Let $B$ be a Banach space.  We will say that an operator $A:
  B\rightarrow B$ is \emph{unitary} if $A$ is invertible and
  $\norm{Ab}=\norm{b}$ for all $b\in B$.  An operator $A: B\rightarrow
  B$ is called \emph{self-adjoint} if the operator $\exp(iA)$ is
  unitary.  In the Hilbert space case this convention coincides with
  the standard definition.
\end{conv}
Let $T_1$, \ldots, $T_n$ be an $n$-tuples of selfadjoint linear
operators on a Banach space $B$. We put for our convenience
$T_0=I$---the identical operator.  It follows from the
Trotter-Daletskii\footnote{The formula is usually attributed to
  Trotter alone.  It is widely unknown that the result appeared
  in~\cite{Daletski60a} also.\label{fn:daletski}}
formula~\cite[Thm.~VIII.31]{SimonReed80} that any linear combination
$\sum_{j=0}^n a_j T_j$ is again a selfadjoint operator.  We will
consider a set of unitary operators
\begin{eqnordisp}[eq:t-set]
  T(a_0,a_1,\ldots,a_n)=\exp\left(i\sum_{j=0}^n a_j T_j \right)
\end{eqnordisp}
parametrized by vectors $(a_0,a_1,\ldots,a_n)\in \Space{R}{n+1}$.
Particularly $T(0,0,\ldots,0)=I$.  A family of their transformations
$\omega(t,z)$, $t\in \Space{R}{}$, $z\in \Space{C}{n}$ is defined by
the rule
\begin{eqnarray}
\omega(t,z)T(a_0,a_1,\ldots,a_n) = T\left(a_0+t +\sum_{j=1}^n(u_j v_j -
\sqrt{2} a_j u_j ), \right. \nonumber\\
\left. \qquad a_1- \sqrt{2}  v_1, \ldots, a_n- \sqrt{2}v_n \right), 
\label{eq:rep-a}
\end{eqnarray}
where $z_j=u_j+iv_j$.  A direct calculation shows that
$\omega(t',z')\omega(t'',z'')=\omega(t'+t''+\frac{1}{2}\Im
(\bar{z}'{z}''),z'+z'')$---this is a non-linear geometric
representation of the Heisenberg group $\Space{H}{n}$.  We could
observe that
\begin{eqnordisp}[]
  T(a_0,a_1,\ldots,a_n)=\omega(a_0, a) T(0,0,\ldots,0) =\omega(a_0, a)
  T_0 =\omega(a_0, a)I,
\end{eqnordisp}
where $a=(ia_1,\ldots,ia_n)$.  Obviously all transformations
$\omega(t,z)$ are isometries if the norm of elements
$T(a_0,a_1,\ldots,a_n)$ is defined as their operator norm.

The representation $\omega$~\eqref{eq:rep-a} is not linear and we
would like to use the procedure outlined in
Remark~\ref{re:linearization}.  We construct the linear space of
operator valued functions $\FSpace{L}{}(\Space{R}{n},\algebra{B})$ for
a $\Space{H}{n}$-homogeneous space $X$ as follows
\begin{eqnordisp}[eq:a-inverse]
  [\oper{T}f](t)=\int_{X} f(x)\,\omega(s(x))\,dx \, T(t), \qquad t \in
  \Space{R}{n},\ s(x)\in \Space{H}{n}.
\end{eqnordisp}
We also extend the representation $\omega$ to
$\FSpace{L}{}(X,\algebra{B})$ in accordance
witn~\eqref{eq:rep-op-ext}.

We will go on with coherent states defined by such a representation.
In the notations of Subsection~\ref{ss:calculus} operators $T_1$,
\ldots, $T_n$ form a set $\mathbf{T}$ defining the representation
$\tau=\omega$ in~\eqref{eq:rep-a} with $T_0(x)=I$ being a vacuum
vector.

\begin{example}\label{ex:weyl} 
  We are ready to demonstrate that the Weyl functional calculus is an
  application of Definition~\ref{de:calc} and Example~\ref{ex:fourier}
  as was announced in~\cite[Remark~4.3]{Kisil95i}. Consider again the
  subgroup $H=\Space{R}{n+1}=\{(t,z)\such \Im(z)=0\}$ and a
  realization of scheme from Subsection~\ref{ss:singular} for this
  subgroup. Then the first paragraph of Example~\ref{ex:fourier} is
  applicable here.
  
  We could easily see that $\omega(t,u_1,\ldots,u_n)I=e^{it+i\sum_1^n
    u_j}I$ and we select the identity operator $I$ times
  $(2\pi)^{-n/2}$ as the vacuum vector $T_0(x)$ of the representation
  $\omega$.  Thereafter the transformation $\oper{T}:
  \FSpace{S}{}(\Space{R}{n}) \rightarrow
  \FSpace{L}{}(\Space{R}{n},\algebra{B})$ \eqref{eq:a-inverse} is
  exactly the inverse wavelet transform for the representation
  $\omega$.  This transformation is defined at least for all $f\in
  \FSpace{S}{}(\Space{R}{n})$.  The space $\FSpace{S}{}(\Space{R}{n})$
  is the image of the wavelet (Fourier) transform~\eqref{eq:fourier}.
  Thus as outlined in Proposition~\ref{pr:calc} we could construct an
  intertwining operator $\fourier{}$ between
  $\sigma$~\eqref{eq:schrodinger} and $\omega$~\eqref{eq:rep-a} from
  the formula \eqref{eq:calc-map} as follow (see
  Remark~\ref{re:zero}):
\begin{eqnarray}
[\Phi_{\mathbf{T}} f] (0) & = & \oper{M}_\omega \oper{W}_\sigma f 
=  (2\pi)^{-n/2} \int_{\Space{R}{n}} \widehat{f}(x)\, \omega(0, x_1,
\ldots, x_n) \, I\,dx \nonumber\\
& = &  (2\pi)^{-n/2} \int_{\Space{R}{n}} \widehat{f}(x)\, e^{i\sum_1^n x_j
T_j} \,dx. \label{eq:weyl-calc}
\end{eqnarray}
This formula is exactly the integral formula for the Weyl functional
calculus~\cite{Anderson69,Nelson70,MTaylor68}. As an example one could
define a function
\begin{eqnordisp}[eq:gauss-op]
  e^{-\sum_1^n T_j^2/2}= (2\pi)^{-n/2} \int_{\Space{R}{n}}
  e^{-\sum_1^n x_j^2/2} e^{i\sum_1^n x_j T_j} \,dx.
\end{eqnordisp}
As we have seen~\eqref{eq:calc-dis} one could formally write the
integral kernel from~\eqref{eq:weyl-calc} as convolution of the
integral kernels for $\oper{W}_\sigma$ and $\oper{M}_\omega$:
\begin{eqnordisp}[eq:weyl-sp]
  \Phi(y,0) = \int_{\Space{R}{n}} e^{-i\sum_1^n y_j x_j} e^{i\sum_1^n
    x_j T_j}\, dx.
\end{eqnordisp}
This expression looks very formal, but it is possible to give it a
precise mathematical meaning as an operator valued distribution.  Such
an approach was explored by \person{Anderson} in~\cite{Anderson69}.
The support of this distribution was defined as the \emph{Weyl joint
  spectrum} for $n$-tuple of non-commuting operators $T_1$,
\ldots,$T_n$ and studied in~\cite{Anderson69}.
\begin{rem}
  As was mentioned in Remark~\ref{re:ax+b} one could construct the
  Fourier transform from representations of $ax+b$ group or the group
  $\Space{M}{n}$ of M\"obius transformations of the upper half plane.
  Analogously one could deduce the Weyl functional calculus as an
  intertwining operator between two representation of this group.  The
  Cauchy kernel $G(x)$~\cite[\S~9]{BraDelSom82} in Clifford analysis
  is the kernel of the Cauchy integral transform
\begin{eqnordisp}[]
  f(y) = \int_{\partial \Omega} G_y (x) \vec{n}(x) f(x)\,d\sigma(x),
  \qquad y\in \Omega,
\end{eqnordisp}
where $\vec{n}(x)$ is the outer unit vector orthogonal to
$\partial\Omega$ and $d\sigma(x)$ is the surface element at a point
$x$.  The Cauchy integral formula (as any wavelet transform)
intertwines two representations acting on $\Omega$ and
$\partial\Omega$ of $\Space{M}{n}$ similarly to the case of complex
analysis~\cite[\S~6]{Kisil97a}.  Thus we could apply here
Theorem~\ref{th:sp-map} on a mapping of spectral distributions.  The
formula~\eqref{eq:sp-map} take the form
\begin{eqnordisp}[]
  \Phi_W f = \int_{\partial \Omega} \Phi_W(G) (x) \vec{n}(x)
  f(x)\,d\sigma(x),
\end{eqnordisp}
where $\Phi_W$ stands for the Weyl functional calculus. This gives
another interpretation for the main result of the
paper~\cite[Thm.~5.4]{JeffMcInt98a}.
\end{rem}
\end{example}

There is no a reason to limit ourself only to the case of subgroup
$H=\Space{R}{n+1}=\{(t,z)\such \Im(z)=0\}$. Thus we proceed with the
next example.
\begin{example} \label{ex:sb-oper}
  In an analogy with Example~\ref{ex:seg-barg} let us consider now the
  wavelet theory associated to the subgroup $H=\Space{R}{1}=\{(t,0)\}$
  and the representation $\omega$.  The first paragraph of
  Example~\ref{ex:seg-barg} depends only on $G=\Space{H}{n}$ and
  $H=\Space{R}{1}$ and thus is applicable in our case.
  
  It is easy to see from formula~\eqref{eq:rep-a} that any operator
  valued function $[\oper{T}f](x)$~\eqref{eq:a-inverse} is an eigen
  vector for $\omega(h)$, $h\in H$.  To be concise with function
  models we select as a vacuum vector the operator $\exp(-\sum_1^n
  T^2_j)$~\eqref{eq:gauss-op}.  Then the
  condition~\eqref{eq:t-vac-def} immediately follows
  from~\eqref{eq:gauss-op}. Thus we could define a functional calculus
  $\Phi: \FSpace{F}{p}(\Space{C}{n})\rightarrow
  \FSpace{L}{}(\Space{C}{n},\algebra{B})$ by the formula (see
  Remark~\ref{re:zero}):
\begin{eqnarray}
[\Phi_{\mathbf{T}}f](0) &=& \int_{\Space{C}{n}} f(z) \omega(0,z) 
\exp(-\sum_1^n
T^2_j/2)\, dz \nonumber\\
 &=& \int_{\Space{C}{n}} f(z) \omega(0,z) (2\pi)^{-\frac{n}{2}}
\int_{\Space{R}{n}} e^{-\sum_1^n x_j^2/2}
e^{i\sum_1^n x_j T_j} \,dx\, dz \nonumber\\
 &=& (2\pi)^{-\frac{n}{2}} \int_{\Space{C}{n}} f(z) \int_{\Space{R}{n}}
e^{-\sum_1^n x_j^2/2}
 \omega(0,z) e^{i\sum_1^n x_j T_j} \,dx\, dz \nonumber\\
 &=& (2\pi)^{-\frac{n}{2}}\! \int_{\Space{R}{3n}}  \exp\sum_{j=1}^n \left( 
- \frac{ x_j^2}{2} +i(v_j-\sqrt[]{2}x_j) u_j +i (x_j - \sqrt[]{2} v_j)
T_j\right)\nonumber\\
&& \qquad \qquad \times f(z)\, dx\,du  \, dv. \label{eq:sb-calc}
\end{eqnarray} 
The last formula could be rewritten for mutually commuting operators
$T_j$ as follows:
\begin{eqnarray}
[\Phi_{\mathbf{T}}f] (0) &=& (2\pi)^{-\frac{n}{2}} \int_{\Space{R}{2n}} 
\int_{\Space{R}{n}} \exp\sum_{j=1}^n \left(- \frac{ 
x_j^2}{2} + ix_j  (T_j-\sqrt[]{2}u_j)\right)\, dx
\nonumber\\
&& \qquad \qquad \times  \exp\sum_{j=1}^n i(v_ju_j- 
\sqrt[]{2} v_jT_j)\,  f(z) \,du  \, dv \nonumber\\
&=& (2\pi)^{-\frac{n}{2}} \int_{\Space{C}{n}} \exp\sum_{j=1}^n \left( 
iv_j(u_j-\sqrt[]{2}T_j) - \frac{ (T_j-\sqrt[]{2}u_j)^2}{2} \right) f(z)\,dz
\nonumber
\end{eqnarray}
where the exponent of operator is defined in the standard sense, e.g.
via the Weyl functional calculus~\eqref{eq:weyl-calc} or the Taylor
expansion.  The last formula is similar to~\eqref{eq:sb-inverse}.
This is very natural for commuting operators as well as that for
non-commuting operators fromula~\eqref{eq:sb-calc} is more
complicated.

The spectral distribution
\begin{eqnordisp}[]
  \Phi_{\mathbf{T}}(z,0)=(2\pi)^{-\frac{n}{2}} \int_{\Space{R}{n}}
  \exp\sum_{j=1}^n \left( - \frac{ x_j^2}{2} +i(v_j-\sqrt[]{2}x_j) u_j
    +i (x_j - \sqrt[]{2} v_j) T_j\right)\,dx
\end{eqnordisp}
derived from~\eqref{eq:sb-calc} contains at least as much information
on operators $T_1$, \ldots, $T_n$ as the Weyl
distribution~\eqref{eq:weyl-sp} and deserves a careful separate
investigation.  We will just mention in conclusion that the
Segal-Bargmann space is an example of the Fock space---space of second
quantization for bosonic fields.  Thus the functional calculus based
on the Segal-Bargmann model sketched here seems to be an appropriate
model for quantized bosonic fields.
\end{example}

\renewcommand{\thesection}{}
\section{Acknowledgments}
The author was supported by the grant INTAS 93--0322--Ext.  I am
grateful to the Prof.~H.~Feichtinger, Prof.~J.~Klauder,
Dr.~D.S.~Kalyuzhny, and Prof.~I.~Segal for their valuable comments on
the subject of the paper.  The history of the Trotter-Daletskii
formula mentioned in Footnote~\ref{fn:daletski} was told me by
Dr.~V.~Kushnirevitch.  I am in debt to Dr.~J.~Cnops, Dr.~M.V.~Kuzmin,
M.Z.~Neyman, and O.P.~Pilipenko who shared with me their philological
knowledge.

\small \bibliographystyle{plain}
\bibliography{abbrevmr,akisil,analyse,acombin,aphysics,arare}

\newcommand{\noopsort}[1]{} \newcommand{\printfirst}[2]{#1}
  \newcommand{\singleletter}[1]{#1} \newcommand{\switchargs}[2]{#2#1}
  \newcommand{\irm}{\textup{I}} \newcommand{\iirm}{\textup{II}}
  \newcommand{\vrm}{\textup{V}} \providecommand{\MathRev}[1]{\textbf{MR}\# #1}
\begin{thebibliography}{10}

\bibitem{AliAntGazMue}
S.~Twareque Ali, J.-P. Antoine, J.-P. Gazeau, and U.A. Mueller.
\newblock Coherent states and their generalizations: {A} mathematical overview.
\newblock {\em Rev. Math. Phys.}, 7(7):1013--1104, 1995.
\newblock Zbl \# 837.43014.

\bibitem{Anderson69}
Robert~F.~V. Anderson.
\newblock The {Weyl} functional calculus.
\newblock {\em J. Funct. Anal.}, 4:240--267, 1969.

\bibitem{Arveson80}
William Arveson.
\newblock {\em The Harmonic Analysis of Automorphisms Groups}.
\newblock American Mathematical Society, Amer. Math. Soc. Summer Institute,
  {\noopsort{}}1980.

\bibitem{AtiyahSchmid80}
Michael Atiyah and Wilfried Schmid.
\newblock A geometric construction of the discrete series for semisimple {Lie}
  group.
\newblock In J.A. Wolf, M.~Cahen, and M.~De Wilde, editors, {\em Harmonic
  Analysis and Representations of Semisimple {Lie} Group}, volume~5 of {\em
  Mathematical Physics and Applied Mathematics}, pages 317--383. D. Reidel
  Publishing Company, Dordrecht, Holland, {\noopsort{}}1980.

\bibitem{Bargmann61}
V.~Bargmann.
\newblock On a {H}ilbert space of analytic functions and an associated integral
  transform. {Part I}.
\newblock {\em Comm. Pure Appl. Math.}, 3:215--228, 1961.

\bibitem{Bargmann67}
V.~Bargmann.
\newblock On a {Hilbert} space of analytic functions and an associated integral
  transform. {II: A} family of related function spaces. {Application} to
  distribution theory.
\newblock {\em Commun. Pure Appl. Math.}, 20:1--101, 1967.
\newblock Zbl \# 149.09601.

\bibitem{Berezin72}
Felix~A. Berezin.
\newblock Covariant and contravariant symbols of operators.
\newblock {\em Izv. Akad. Nauk SSSR Ser. Mat.}, 6:1117--1151, 1972.
\newblock Reprinted in~\cite[pp.~228--261]{Berezin88}.

\bibitem{Berezin74}
Felix~A. Berezin.
\newblock Quantization.
\newblock {\em Math. USSR-Izv.}, 8:1109--1165, 1974.

\bibitem{Berezin75}
Felix~A. Berezin.
\newblock Quantization in complex symmetric spaces.
\newblock {\em Math. USSR-Izv.}, 9:341--379, 1975.

\bibitem{Berezin88}
Felix~A. Berezin.
\newblock {\em Method of Second Quantization}.
\newblock ``Nauka'', Moscow, {\noopsort{}}1988.

\bibitem{BernTayl94}
David Bernier and Keith~F. Taylor.
\newblock Wavelets from square-integrable representations.
\newblock {\em SIAM J. Math. Anal.}, 27(2):594--608, 1996.
\newblock \MathRev{97h:22004}.

\bibitem{BraDelSom82}
F.~Brackx, R.~Delanghe, and F.~Sommen.
\newblock {\em Clifford Analysis}, volume~76 of {\em Research Notes in
  Mathematics}.
\newblock Pitman Advanced Publishing Program, Boston, 1982.

\bibitem{BratJorg97a}
Ola Bratteli and Palle E.~T. Jorgensen.
\newblock Isometries, shifts, {C}untz algebras and multiresolution wavelet
  analysis of scale ${N}$.
\newblock {\em Integral Equations Operator Theory}, 28(4):382--443, 1997.
\newblock \eprint{http://xxx.lanl.gov/abs/funct-an/9612003}{funct-an/9612003}.

\bibitem{Cnops94a}
Jan Cnops.
\newblock {\em {Hurwitz} Pairs and Applications of {M\"obius} Transformations}.
\newblock {Habilitation} dissertation, Universiteit Gent, Faculteit van de
  Wetenschappen, 1994.
\newblock \eprint{ftp://cage.rug.ac.be/pub/clifford/jc9401.tex}
  {ftp://cage.rug.ac.be/pub/clifford/jc9401.tex}.

\bibitem{CnopsKisil97a}
Jan Cnops and Vladimir~V. Kisil.
\newblock Monogenic functions and representations of nilpotent {Lie} groups in
  quantum mechanics.
\newblock {\em Mathematical Methods in the Applied Sciences}, 22(4):353--373,
  1998.
\newblock \eprint{http://xxx.lanl.gov/abs/math/9806150/}{math/9806150}.

\bibitem{Coburn94a}
Lewis~A. Coburn.
\newblock {Berezin-Toeplitz} quantization.
\newblock In {\em Algebraic Mettods in Operator Theory}, pages 101--108.
  Birkh\"auser Verlag, New York, 1994.

\bibitem{Connes73}
Alain Connes.
\newblock Une classification des facteurs de type{ \textup{III}}.
\newblock {\em Ann. Sci. {\'E}cole Norm. Sup. (4)}, 6:133--252, 1973.

\bibitem{Cordes79}
H.~O. Cordes.
\newblock {\em Pseudodifferential Operators --- An Abstract Theory}, volume 756
  of {\em Lect. Notes Math.}
\newblock Springer-Verlag, Berlin, 1979.

\bibitem{CotlSad90b}
Mischa Cotlar and Cora Sadosky.
\newblock Toeplitz and {Hankel} forms related to unitary representations of the
  symplectic plane.
\newblock {\em Colloq. Math.}, 50/51:693--708, 1990.

\bibitem{Daletski60a}
Yu.L. Daletski.
\newblock On representation of solutions of operator equations in a form of
  functional integrals.
\newblock {\em Dokl. Akad. Nauk SSSR}, 134(5):1013--1016, 1960.

\bibitem{Daubechies92}
Ingrid Daubechies.
\newblock {\em Ten Lectures on Wavelets}, volume~61 of {\em CBMS-NSF Regional
  Conference Series in Applied Mathematics}.
\newblock Society for Industrial and Applied Mathematics (SIAM), Philadelphia,
  PA, {\noopsort{}}1992.

\bibitem{DelSomSou92}
Richard Delanghe, Frank Sommen, and Vladimir Sou\v{c}ek.
\newblock {\em Clifford Algebra and Spinor-Valued Functions}.
\newblock Kluwer Academic Publishers, Dordrecht, {\noopsort{}}1992.

\bibitem{FeichGroech89a}
{Feichtinger, Hans G. and Groechenig, K.H.}
\newblock Banach spaces related to integrable group representations and their
  atomic decompositions, {I}.
\newblock {\em J. Funct. Anal.}, 86(2):307--340, 1989.
\newblock Zbl \# 691.46011.

\bibitem{FeichGroech89b}
{Feichtinger, Hans G. and Groechenig, K.H.}
\newblock Banach spaces related to integrable group representations and their
  atomic decompositions. {Part II}.
\newblock {\em Monatsh. Math.}, 108:129--148, 1989.

\bibitem{Folland89}
Gerald~B. Folland.
\newblock {\em Harmonic Analysis in Phase Space}.
\newblock Princeton University Press, Princeton, New Jersey, {\noopsort{}}1989.

\bibitem{Heidegger61a}
Martin Heidegger.
\newblock {\em \~{V}om Wesen der Wahrheit}.
\newblock Frankfurt am Main, {\noopsort{}}1961.

\bibitem{HeilWaln89}
Christopher~E. Heil and David~F. Walnut.
\newblock Continuous and discrete wavelet transforms.
\newblock {\em SIAM Rev.}, 31(4):628--666, 1989.

\bibitem{Howe80a}
Roger Howe.
\newblock On the role of the {Heisenberg} group in harmonic analysis.
\newblock {\em Bull. Amer. Math. Soc. (N.S.)}, 3(2):821--843, 1980.

\bibitem{Howe80b}
Roger Howe.
\newblock Quantum mechanics and partial differential equations.
\newblock {\em J. Funct. Anal.}, 38:188--254, 1980.

\bibitem{JeffMcInt98a}
Brian Jefferies and Alan McIntosh.
\newblock The {Weyl} calculus and {Clifford} analysis.
\newblock {\em Bull. Austral. Math. Soc.}, 57(2):329--341, 1998.

\bibitem{JorgWer94a}
Palle E.~T. J{\o}rgensen and R.~F. Werner.
\newblock Coherent states of the $q$-canonical commutation relation.
\newblock {\em Comm. Math. Phys.}, 164:455--471, 1994.

\bibitem{KantAkil84}
L.~V. Kantorovich and G.~P. Akilov.
\newblock {\em Functional analysis}.
\newblock Pergamon Press, Oxford, second edition, 1982.
\newblock Translated from the Russian by Howard L. Silcock.

\bibitem{KirGvi82}
A.~A. Kirillov and A.~D. Gvishiani.
\newblock {\em Theorems and Problems in Functional Analysis}.
\newblock Problem Books in Mathematics. Springer-Verlag, New York,
  {\noopsort{}1982}.

\bibitem{Kirillov76}
Alexander~A. Kirillov.
\newblock {\em Elements of the Theory of Representations}, volume 220 of {\em A
  Series of Comprehensive Studies in Mathematics}.
\newblock Springer-Verlag, New York, 1976.

\bibitem{Kisil95a}
Vladimir~V. Kisil.
\newblock Integral representation and coherent states.
\newblock {\em Bull. Belg. Math. Soc. Simon Stevin}, 2(5):529--540, 1995.
\newblock \MathRev{97b:22012}.

\bibitem{Kisil95i}
Vladimir~V. Kisil.
\newblock M\"obius transformations and monogenic functional calculus.
\newblock {\em \href{http://www.ams.org/era/}{Electron. Res. Announc. Amer.
  Math. Soc.}},
  2(1):\href{http://www.ams.org/jourcgi/jour--pbprocess?fn=110&arg1=S1079--676%
2--96--00004--2&u=/era/1996--02--01/}{26--33}, 1996.
\newblock (electronic) \MathRev{98a:47018}.

\bibitem{Kisil97b}
Vladimir~V. Kisil.
\newblock The umbral calculus: a model from convoloids.
\newblock page~25, 1997.
\newblock \eprint{http://xxx.lanl.gov/abs/funct-an/9704001/}{funct-an/9704001}.
  Submitted to \textit{Z. Anal. Anwendungen}.

\bibitem{Kisil97c}
Vladimir~V. Kisil.
\newblock Analysis in{ $\Space{R}{1,1}$} or the principal function theory.
\newblock {\em Complex Variables Theory Appl.}, page~25, 1998.
\newblock (To
  appear)\eprint{http://xxx.lanl.gov/abs/funct-an/9712003/}{funct-an/9712003}.

\bibitem{Kisil95e}
Vladimir~V. Kisil.
\newblock Spectrum of operators, functional calculi and group representations.
\newblock 1999.
\newblock (In preparation).

\bibitem{Kisil97a}
Vladimir~V. Kisil.
\newblock Two approaches to non-commutative geometry.
\newblock In H.~Begehr, O.~Celebi, and W.~Tutschke, editors, {\em Complex
  Methods for Partial Differential Equations}, chapter~14, pages 219--248.
  Kluwer Academic Publishers, Netherlands, 1999.
\newblock \eprint{http://xxx.lanl.gov/abs/funct-an/9703001/}{funct-an/9703001}.

\bibitem{Klauder94b}
John~R. Klauder.
\newblock Coherent states and coordinate-free quantization.
\newblock {\em International Journal of Theoretical Physics}, 33(3):509--522,
  1994.

\bibitem{Klauder94a}
John~R. Klauder.
\newblock Universal propagator for group-related coherent states.
\newblock {\em Acta Physica Polonica A}, 85(4):655--666, 1994.

\bibitem{Klauder95a}
John~R. Klauder.
\newblock Quantization without quantization.
\newblock {\em Annals of Physics}, 237(1):147--160, 1995.

\bibitem{KlauderMcKenna67}
John~R. Klauder and Jim McKenna.
\newblock Wavelets in the space of nontempered distributions.
\newblock Unpublished, 1967.

\bibitem{KlaSkag85}
John~R. Klauder and Bo-Sture Skagerstam, editors.
\newblock {\em Coherent States. Applications in physics and mathematical
  physics.}
\newblock World Scientific Publishing Co., Singapur, {\noopsort{}}1985.

\bibitem{KnappWallach76}
A.W. Knapp and N.R. Wallach.
\newblock Szeg\"o kernels associated with discrete series.
\newblock {\em Invent. Math.}, 34(2):163--200, 1976.

\bibitem{Krein48a}
M.~G. Kre{\u\i}n.
\newblock On {H}ermitian operators with directed functionals.
\newblock {\em Akad. Nauk Ukrain. RSR. Zbirnik Prac\cprime\ Inst. Mat.},
  1948(10):83--106, 1948.
\newblock \MathRev{14:56c}, reprinted in~\cite{KreinII}.

\bibitem{KreinII}
M.~G. Kre{\u\i}n.
\newblock {\em {\cyr {I}zbrannye Trudy}. {II}}.
\newblock Akad. Nauk Ukrainy Inst. Mat., Kiev, 1997.
\newblock \MathRev{96m:01030}.

\bibitem{KreinRutman48}
M.~G. Kre{\u\i}n and M.~A. Rutman.
\newblock Linear operators leaving invariant a cone in a {B}anach space.
\newblock {\em Uspehi Matem. Nauk (N. S.)}, 3(1(23)):3--95, 1948.
\newblock \MathRev{10:256c}.

\bibitem{Rota95}
Joseph~P.S. Kung, editor.
\newblock {\em Gian-Carlo Rota on Combinatorics: Introductory Papers and
  Commentaries}, volume~1 of {\em Contemporary Mathematicians}.
\newblock Birkh\"auser Verlag, Boston, 1995.

\bibitem{FoundationIII}
Ronald Mullin and Gian-Carlo Rota.
\newblock On the foundation of combinatorial theory ({III}): Theory of binomial
  enumeration.
\newblock In B.Harris, editor, {\em Graph Theory and Its Applications}, pages
  167--213. Academic Press, Inc., New York, 1970.
\newblock Reprinted in~\cite[pp.~118--147]{Rota95}.

\bibitem{Nelson70}
E.~Nelson.
\newblock A functional calculus for non-commuting operators.
\newblock In F.E. Browder, editor, {\em Functional analysis and related fields,
  Proceedings of a conference in honour of Professor Marshal Stone, Univ. of
  Chicago, May 1968}, pages 172--187. Springer-Verlag, Berlin, Heidelberg, New
  York, {\noopsort{}}1970.

\bibitem{Perelomov86}
A.~M. Perelomov.
\newblock {\em Generalized Coherent States and Their Applications}.
\newblock Springer-Verlag, Berlin, {\noopsort{}}1986.

\bibitem{SimonReed80}
Michael Reed and Barry Simon.
\newblock {\em Functional Analysis}, volume~1 of {\em Methods of Modern
  Mathematical Physics}.
\newblock Academic Press, Orlando, second edition, {\noopsort{1974}}1980.

\bibitem{RomRota78}
S.~Roman and Gian-Carlo Rota.
\newblock The umbral calculus.
\newblock {\em Adv. in Math.}, 27:95--188, 1978.

\bibitem{Rota64a}
Gian-Carlo Rota.
\newblock The number of partitions of a set.
\newblock {\em Amer. Math. Monthly}, 71(5):498--504, May 1964.
\newblock Reprinted in~\cite[pp.~1--6]{Rota75} and~\cite[pp.~111--117]{Rota95}.

\bibitem{Rota75}
Gian-Carlo Rota.
\newblock {\em Finite Operator Calculus}.
\newblock Academic Press, Inc., New York, 1975.

\bibitem{KahOdlRota73}
Gian-Carlo Rota, David Kahaner, and Andrew Odlyzko.
\newblock Finite operator calculus.
\newblock {\em J. Math. Anal. Appl.}, 42(3):685--760, June 1973.
\newblock Reprinted in~\cite[pp.~7--82]{Rota75}.

\bibitem{Segal60}
Irving~E. Segal.
\newblock {\em Mathematical Problems of Relativistic Physics}, volume~II of
  {\em Proceedings of the Summer Seminar (Boulder, Colorado, 1960)}.
\newblock American Mathematical Society, Providence, R.I., 1963.

\bibitem{Segal90}
Irving~E. Segal.
\newblock Algebraic quantization and stability.
\newblock In William~B. Arveson and Ronald~G. Douglas, editors, {\em Operator
  Theory: Operator Algebras and Applications}, number~51 in Proceedings of
  Symposia in Mathematics, pages 503--518. American Mathematical Society,
  Providence, R.I., 1990.

\bibitem{Segal90a}
Irving~E. Segal.
\newblock The mathematical implications of fundamental physical principles.
\newblock In {\em The Legacy of {John v. Neumann} ({Providence, 1990})},
  volume~50 of {\em Proc. Sympos. Pure Math.}, pages 151--178. American
  Mathematical Society, Providence, R.I., 1990.

\bibitem{Segal94}
Irving~E. Segal.
\newblock {$C^*$}-algebras and quantization.
\newblock In Robert Doran, editor, {\em {$C^*$}-Algebras: 1943--1993}, number
  167 in Contemporary Mathematics, pages 55--65. American Mathematical Society,
  Providence, R.I., 1994.

\bibitem{Segal96a}
Irving~E. Segal.
\newblock Rigorous covariant form of the correspondence principle.
\newblock In William Arveson et~al., editors, {\em Quantization, Nonlinear
  Partial Differential Equations, and Operator Algebra ({Cambridge, MA},
  1994)}, volume~59 of {\em Proc. Sympos. Pure Math.}, pages 175--202. American
  Mathematical Society, Providence, R.I., 1996.

\bibitem{Shubin87}
Mikhail~A. Shubin.
\newblock {\em Pseudodifferential Operators and Spectral Theory}.
\newblock Springer-Verlag, Berlin, 1987.

\bibitem{Sommen88}
Frank Sommen.
\newblock Plane wave decompositions of monogenic functions.
\newblock {\em Annales Pol. Math.}, 49:101--114, 1988.

\bibitem{Takesaki83}
Masamichi Takesaki.
\newblock {\em Structure of Factors and Automorphism Groups}, volume~51 of {\em
  Regional Conference Series in Mathematics}.
\newblock American Mathematical Society, Providence, R.I., {\noopsort{}}1983.

\bibitem{MTaylor68}
Michael~E. Taylor.
\newblock Functions of several self-adjoint operators.
\newblock {\em J. Funct. Anal.}, 19:91--98, 1968.

\bibitem{MTaylor81}
Michael~E. Taylor.
\newblock {\em Pseudodifferential Operators}, volume~34 of {\em Princeton
  Mathematical Series}.
\newblock Princeton University Press, Princeton, New Jersey, 1981.

\bibitem{MTaylor86}
Michael~E. Taylor.
\newblock {\em Noncommutative Harmonic Analysis}, volume~22 of {\em Math. Surv.
  and Monographs}.
\newblock American Mathematical Society, Providence, R.I., {\noopsort{}}1986.

\end{thebibliography}
\end{document}